\documentclass[aos,seceqn,nameyear,dvips]{arximspdf}
\usepackage{dcolumn}
\usepackage{graphicx}

\doi{10.1214/10-AOS819}
\volume{38}
\issue{6}
\pubyear{2010}
\firstpage{3352}
\lastpage{3386}

\makeatletter

\newcolumntype{d}[1]{D{.}{.}{#1}}

\newcommand{\eqref}[1]{(\ref{#1})}

\newtheorem{theorem}{Theorem}
\newtheorem{lemma}{Lemma}
\newtheorem{proposition}{Proposition}
\newproclaim{remark}{Remark}

\newcommand{\vectornorm}[1]{\|#1\|}
\newcommand{\uniform}[1]{\|#1\|_\calB}
\newcommand{\hilbert}[1]{\|#1\|_\calS}

\newcommand{\htheta}{\hat{\theta}}

\newcommand{\norm}[1]{\|#1\|}

\newcommand{\inner}[2]{\langle #1, #2 \rangle}

\newcommand{\kernel}[1]{\operatorname{Ker}(#1)}

\newcommand{\image}[1]{\operatorname{Im}(#1)}

\newcommand{\cov}{\operatorname{Cov}}

\newcommand{\diag}{\operatorname{diag}}

\newcommand{\var}{\operatorname{Var}}

\newcommand{\ve}{{\varepsilon}}
\newcommand{\bA}{{\mathbf A}}
\newcommand{\bB}{{\mathbf B}}
\newcommand{\bK}{{\mathbf K}}
\newcommand{\bW}{{\mathbf W}}
\newcommand{\bY}{{\mathbf Y}}
\newcommand{\be}{{\mathbf e}}

\newcommand{\bbeta}{\bolds{\beta}}
\newcommand{\bfeta}{{\bolds{\eta}}}
\newcommand{\bvarphi}{\bolds{\varphi}}

\newcommand{\bSigma}{\bolds{\Sigma}}
\newcommand{\bgamma}{\bolds{\gamma}}
\newcommand{\brho}{\bolds{\rho}}
\newcommand{\bxi}{\bolds{\xi}}

\newcommand{\calI}{{\mathcal I}}
\newcommand{\calL}{{\mathcal L}}
\newcommand{\calX}{{\mathcal X}}
\newcommand{\calY}{{\mathcal Y}}
\newcommand{\calV}{{\mathcal V}}
\newcommand{\calZ}{{\mathcal Z}}
\newcommand{\calF}{{\mathcal F}}
\newcommand{\calH}{{\mathcal H}}

\newcommand{\calM}{{\mathcal M}}
\newcommand{\calN}{{\mathcal N}}
\newcommand{\calP}{{\mathcal P}}
\newcommand{\calS}{{\mathcal S}}
\newcommand{\calB}{{\mathcal B}}

\makeatother

\begin{document}
\begin{frontmatter}

\title{Identifying the finite dimensionality of~curve~time~series\thanksref{T1}}
\runtitle{Curve time series}

\thankstext{T1}{Supported in part by the Engineering and Physical Sciences Research Council of the UK.}

\begin{aug}
\author[A]{\fnms{Neil} \snm{Bathia}\ead[label=e1]{n.bathia@lse.ac.uk}},
\author[A]{\fnms{Qiwei} \snm{Yao}\corref{}\ead[label=e2]{q.yao@lse.ac.uk}} and
\author[B]{\fnms{Flavio} \snm{Ziegelmann}\ead[label=e3]{flavioaz@mat.ufrgs.br}}
\runauthor{N. Bathia, Q. Yao and F. Ziegelmann}
\affiliation{London School of Economics, London School of Economics and
Federal~University~of~Rio Grande do Sul}
\address[A]{N. Bathia\\
Q. Yao\\
Department of Statistics\\
London School of Economics\\
London, WC2A 2AE\\
United Kingdom\\
\printead{e1}\\
\phantom{E-mail: }\printead*{e2}}
\address[B]{F. Ziegelmann\\
Department of Statistics\\
Federal University\\
\quad of Rio Grande do Sul\\
91509-900 Porto Alegre\\
Brazil\\
\printead{e3}}
\end{aug}

\received{\smonth{10} \syear{2009}}
\revised{\smonth{2} \syear{2010}}

%
\begin{abstract}
The curve time series framework provides a convenient vehicle to
accommodate some nonstationary features into a stationary setup. We
propose a new method to identify the dimensionality of curve time
series based on the dynamical dependence across different curves. The
practical implementation of our method boils down to an eigenanalysis
of a finite-dimensional matrix. Furthermore, the determination of the
dimensionality is equivalent to the identification of the nonzero
eigenvalues of the matrix, which we carry out in terms of some
bootstrap tests. Asymptotic properties of the proposed method are
investigated. In particular, our estimators for zero-eigenvalues enjoy
the fast convergence rate $n$ while the estimators for nonzero
eigenvalues converge at the standard $\sqrt n$-rate. The proposed
methodology is illustrated with both simulated and real data sets.
\end{abstract}

%

\end{frontmatter}

%
\section{Introduction}\label{section:2.1}
%

A curve time series may consist of, for example, annual weather
record charts, annual production charts or daily volatility curves
(from morning to evening). In these examples, the curves are segments
of a single long time series. One advantage to view them as a curve
series is to accommodate some nonstationary features (such as
seasonal cycles or diurnal volatility patterns) into a stationary
framework in a Hilbert space.
There are other types of curve series
that cannot be pieced together into a single long time series; for
example, daily mean-variance efficient frontiers of portfolios,
yield curves and intraday asset return distributions. See also an example
of daily return density curves in Section~\ref{section:2.4.2}. The goal of
this paper is to identify the finite dimensionality of curve time
series in the sense that the serial dependence across different
curves is driven by a finite number of scalar components. Therefore,
the problem of modeling curve dynamics is reduced to that of
modeling a finite-dimensional vector time series.

Throughout this paper, we assume that the observed curve time series,
which we denote by $Y_1(\cdot),\ldots, Y_n(\cdot)$, are defined on
a compact interval $\calI$ and are subject to errors in the sense
that
%
\begin{equation} \label{a1}
Y_t (u) = X_t(u) + \ve_t(u),\qquad u\in\calI,
\end{equation}
where $ X_t(\cdot) $ is the curve process of interest. The existence
of the noise term $\ve_t(\cdot)$ reflects the fact that curves
$X_t(\cdot)$ are seldom perfectly observed. They are often only
recorded on discrete grids and are subject to both experimental
error and numerical rounding. These noisy discrete data are smoothed
to yield ``observed'' curves $Y_t(\cdot)$. Note that both $X_t(\cdot)$
and $\ve_t(\cdot)$ are unobservable.

We assume that $\ve_t(\cdot)$ is a white noise sequence in the sense
that $E\{\ve_t(u)\} = 0$ for all $t$ and $\cov\{\ve_t(u), \ve_s(v)\}
= 0$ for any $u,v \in\calI$ provided $t \ne s$. This is
guaranteed since we may include all the dynamic elements of
$Y_t(\cdot)$ into $X_t(\cdot)$. Likewise, we may also assume that no
parts of $X_t(\cdot)$ are white noise since these parts should be
absorbed into $\ve_t(\cdot)$. We also assume that
%
\begin{equation} \label{a2}
\int_\calI E \{ X_t(u)^2 + \ve_t(u)^2 \} \,du < \infty,
\end{equation}
and both
%
\begin{equation} \label{a21}
\mu(u) \equiv E\{ X_t(u) \},\qquad M_{k}(u,v) \equiv\cov\{
X_t(u), X_{t+k}(v) \}
\end{equation}
do not depend on $t$. Furthermore, we assume that $X_t(\cdot)$ and
$\ve_{t+k}(\cdot)$ are uncorrelated for all integer $k$.
Under condition (\ref{a2}),
$X_t(\cdot)$ admits the Karhunen--Lo\'eve expansion
%
\begin{equation} \label{a22}
X_t(u) - \mu(u) = \sum_{j=1}^\infty\xi_{tj}
\varphi_j(u),
\end{equation}
where $\xi_{tj} = \int_\calI \{ X_t(u)
- \mu(u)\} \varphi_j(u) \,du$ with $\{\xi_{tj}, j\ge1\}$ being a
sequence of scalar random variables with $E(\xi_{tj}) = 0$,
$\var(\xi_{tj}) = \lambda_j$ and $\cov(\xi_{ti}, \xi_{tj}) = 0$ if
$i \neq j$. We rank $\{\xi_{tj}, j\ge1\}$ such that $\lambda_j$
is monotonically decreasing as $j$ increases.

We say that $X_t(\cdot)$ is $d$-dimensional if $\lambda_d\ne0$ and
$\lambda_{d+1} = 0 $, where $d\ge1$ is a finite integer; see
\citet{hall2006}. The primary goal of this paper is to identify $d$
and to estimate the dynamic space $\calM$ spanned by the
(deterministic) eigenfunctions $\varphi_1(\cdot),\ldots,
\varphi_d(\cdot)$.

\citet{hall2006} tackle this problem under the assumption that the
curves $Y_1(\cdot),\ldots, Y_n(\cdot)$ are independent. Then the
problem is insoluble in the sense that one cannot separate
$X_t(\cdot)$ from $\ve_t(\cdot)$ in (\ref{a1}). This difficulty was
resolved in \citet{hall2006} under a ``low noise'' setting which
assumes that the noise $\ve_t(\cdot)$ goes to 0 as the sample size
goes to infinity. Our approach is different and it does not require
the ``low noise'' condition, since we identify $d$ and $\calM$ in
terms of the serial dependence of the curves. Our method relies on a
simple fact that $M_k(u,v) = \cov\{ Y_t(u), Y_{t+k}(v) \}$ for any
$k\ne0$, which automatically filters out the noise $\ve_t(\cdot)$;
see (\ref{a21}). In this sense, the existence of dynamic dependence
across different curves makes the problem tractable.

Dimension reduction plays an important role in functional data
analysis. The most frequently used method is the functional principal
component analysis in the form of applying the Karhunen--Lo\'eve
decomposition directly to the observed curves. The literature in
this field is vast and includes
\citet{besse1986}, \citet{dauxois1982}, \citet{ramsay1991},
\citet{rice1991} and \citet{ramsay2005}.
In spite of the methodological advancements
with independent observations, the work
on functional time series has been of a more theoretical nature; see,
for example, \citet{Bosq2000}. The available inference methods focus mostly
on nonparametric estimation for some characteristics of functional
series [Part IV of \citet{fv2006}]. As far as we are aware, the
work presented here represents the first attempt on the
dimension reduction based
on dynamic dependence, which is radically different from the existing
methods.
Heuristically, our
approach differs from functional principal components analysis in one
fundamental manner; in principal component analysis the objective is to
find the linear combinations of the data which maximize variance. In
contrast, we seek for the linear combinations of the data which
represent the serial dependence in the data.
Although we confine ourselves to square
integrable curve series in this paper, the methodology may be
extended to a more general functional framework including, for
example, a surface series which is particularly important for
environmental study; see, for example, \citet{guillas2008}. A follow-up
study in this
direction will be reported elsewhere.

The rest of the paper is organized as follows. Section
\ref{section:2.2} introduces the proposed new methodology
for identifying the finite-dimensional dynamic structure.
Although the Karhunen--Lo\'eve decomposition (\ref{a22}) serves as a
starting point, we do not seek for such a decomposition explicitly.
Instead the eigenanalysis is performed on a positive-definite operator
defined based on the autocovariance function of the curve process.
Furthermore, computationally our method boils
down to an eigenanalysis of a finite matrix thus requiring no
computing of eigenfunctions in a functional space directly. The
relevant theoretical results are presented in Section
\ref{section:2.3}. As our estimation for the eigenvalues
are essentially quadratic, the convergence rate
of the estimators for the zero-eigenvalues is $n$ while
that for the nonzero eigenvalues is standard $\sqrt{n}$.
Numerical illustration using both simulated and
real datasets is provided in Section~\ref{section:2.4}. Given the
nature of the subject concerned, it is inevitable to make use of
some operator theory in a Hilbert space. We collect some relevant
facts in Appendix~\ref{appA}. We relegate all the technical proofs to
Appendix~\ref{appB}.

%
\section{Methodology}\label{section:2.2}
%

\subsection{Characterize $d$ and $\calM$ via serial dependence}
\label{section:2.2.1}

Let $\calL_2(\calI)$ denote the Hilbert space consisting of all the
square integrable curves defined on $\calI$ equipped with the inner
product
%
\begin{equation}\label{a36}
\inner{f} {g} = \int_\calI f(u) g(u) \,du,\qquad f,g\in \calL_2(\calI).
\end{equation}
Now $M_k$ defined in (\ref{a21}) may be viewed as the kernel of a
linear operator acting on $\calL_2(\calI)$, that is, for any $g\in
\calL_2(\calI)$, $M_k$ maps $g(u)$ to $ \check g(u) \equiv
\int_\calI M_k(u,v) g(v) \,dv$. For notational economy, we will use
$M_k$ to denote both the kernel and the operator. Appendix~\ref{appA} lists
some relevant facts about operators in Hilbert spaces.

For $M_0$ defined in (\ref{a21}), we have a spectral decomposition
of the form
%
\begin{equation}
\label{b1} M_0(u,v) = \sum_{j=1}^\infty\lambda_j \varphi_j(u)
\varphi_j(v),\qquad u,v \in\calI,
\end{equation}
where $\lambda_1\ge\lambda_2 \ge\cdots\ge0$ are the eigenvalues and
$\varphi_1, \varphi_2, \ldots$ are the corresponding orthonormal
eigenfunctions (i.e., $ \inner{\varphi_i}{\varphi_j}
=1 $ for $i=j$, and 0 otherwise). Hence,
\[
\int_\calI M_0(u,v) \varphi_j(v) \,dv = \lambda_j \varphi_j(u),\qquad
j \geq1.
\]
Furthermore, the random curves $X_t(\cdot)$ admit the representation
(\ref{a22}). We assume in this paper that $X_t(\cdot)$ is
$d$-dimensional (i.e., $\lambda_{d+1}=0$). Therefore,
%
\begin{equation}\label{a4}
M_0(u,v) = \sum_{j=1}^d \lambda_j \varphi_j(u) \varphi_j(v),\qquad
X_t(u) = \mu(u) + \sum_{j=1}^d \xi_{tj} \varphi_j(u).
\end{equation}
It follows from (\ref{a1}) that
%
\begin{equation}
\label{a5}
Y_t(u) = \mu(u) + \sum_{j=1}^d \xi_{tj} \varphi_j(u) + \ve_t(u).
\end{equation}
Thus, the serial dependence of $Y_t(\cdot)$ is determined entirely
by that of the $d$-vector process $\bxi_t \equiv(\xi_{t1},\ldots,
\xi_{td})^\prime$ since $\ve_t(\cdot)$ is white noise. By the virtue
of the Karhunen--Lo\'eve decomposition, $E\bxi_t = 0$ and Var$(\bxi_t)
= \diag(\lambda_1,\ldots, \lambda_d)$.

For some prescribed integer $p$, let
%
\begin{equation}
\label{kernel} \widehat{M}_k(u,v) = {1\over n-p} \sum_{j=1}^{n-p}
\{Y_j(u) - \bar{Y}(u) \} \{Y_{j+k}(v) -
\bar{Y}(v) \},
\end{equation}
where $\bar Y(\cdot) = n^{-1} \sum_{1\le j \le n} Y_j (\cdot)$ and
$k=1,\ldots,p$. The reason for truncating the sums
in (\ref{kernel}) at $n-p$ as opposed to $n-k$ is to ensure a
duality operation which simplifies the computation for
eigenfunctions; see Remark~\ref{remark2} at the end of
Section~\ref{section:2.2.2.2}. The conventional approach to estimate $d$ and
$\calM= \operatorname{span}\{ \varphi_1(\cdot), \ldots,
\varphi_d(\cdot)\}$ is to perform an eigenanalysis on $\widehat M_0$ and
let $\widehat d$ be the number
of nonzero eigenvalues and $\widehat\calM$ be spanned by the $\widehat d$
corresponding eigenfunctions; see, for example, \citet{ramsay2005}
and references therein. However, this approach suffers from
complications due to fact that $\widehat M_0$ is not a consistent
estimator for $M_0$, as Cov$\{ Y_t(u), Y_t(v) \} = M_0(u,v) + \cov\{
\ve_t(u), \ve_t(v) \}$. Therefore, $\widehat M_0$ needs to be adjusted to
remove the part due to $\ve_t(\cdot)$ before the eigenanalysis may
be performed. Unfortunately, this is a nontrivial matter since
both $X_t(\cdot)$ and $\ve_t(\cdot)$ are unobservable. An alternative
is to
let the variance of $\ve_t(\cdot)$ decay to 0 as the sample size $n$
goes to infinity; see \citet{hall2006}.

We adopt a different approach based on the fact that Cov$\{ Y_t(u),
Y_{t+k}(v) \} = M_k(u,v)$ for any $k\ne0$, which ensures that $\widehat
M_k$ is a legitimate estimator for $M_k$; see (\ref{a21}) and
(\ref{kernel}).

Let $\bSigma_k = E(\bxi_t \bxi_{t+k}') \equiv( \sigma_{ij}^{(k)})$
be the autocovariance matrix of $\bxi_t$ at lag $k$.
It is easy to see from (\ref{a21}) and (\ref{a4}) that
$M_k(u, v) =\sum_{i,j=1}^d \sigma_{ij}^{(k)}\varphi_i(u)
\varphi_j(v)$. Define a nonnegative operator
%
\begin{equation} \label{a34}
N_k(u,v) = \int_\calI M_k(u,z) M_k(v,z) \,dz =
\sum_{i,j=1}^d w_{ij}^{(k)} \varphi_i(u) \varphi_j(v),
\end{equation}
where $ \bW_k = (w_{ij}^{(k)} ) = \bSigma_k \bSigma_k^\prime$ is a
nonnegative definite matrix. Then it holds for any integer $k$ that
%
\begin{equation} \label{a30}
\int_\calI N_k(u, v) \zeta(v) \,dv = 0\qquad \mbox{for any }
\zeta(\cdot) \in\calM^\bot,
\end{equation}
where $\calM^\bot$ denotes the orthogonal complement of $\calM$ in
$\calL_2(\calI)$. Note (\ref{a30}) also holds if we replace $N_k$ by
the operator
%
\begin{equation} \label{a31}
K(u,v) = \sum_{k=1}^p N_k(u, v),
\end{equation}
which is also a nonnegative operator on $\calL_2(\calI)$.
\begin{proposition}\label{prop:1}
Let the matrix $\bSigma_{k_0}$ be full-ranked for some $ k_0\ge1$.
Then the assertions below hold.

\begin{longlist}
\item The operator $N_{k_0}$ has exactly $d$ nonzero eigenvalues, and
$\calM$ is the linear space spanned by the corresponding $d$
eigenfunctions.

\item For $p\ge k_0$, \textup{(i)} also holds for the operator $K$.
\end{longlist}
\end{proposition}
\begin{remark}\label{remark1}
(i) The condition that rank$(\bSigma_k)
=d$ for some $k\ge1$ is implied by the assumption that $X_t(\cdot)$
is $d$-dimensional. In the case where rank$(\bSigma_k) <d$ for all
$k$, the component with no serial correlations in $X_t(\cdot)$
should be absorbed into white noise $\ve_t(\cdot)$; see similar
arguments on modeling vector time series in \citet{pena1987} and
\citet{pan2008}.

\mbox{\phantom{i}}(ii) The introduction of the operator $K$ in (\ref{a31}) is to pull
together the information at different lags. Using single $N_k$ may
lead to spurious choices of $\widehat d$.

(iii) Note that $\int_\calI K(u,v) \zeta(v)\,dv =0$ if and only if
$\int_\calI N_k(u,v) \zeta(v)\,dv =0$ for all $1\le k\le p$. However,
we cannot use $M_k$ directly in defining $K$ since it does not
necessarily hold that $\int_{\calI} \sum_{1 \leq k \leq p} M_k(u,v)
g(v) \neq0$ for all $g \in\calM$. This is due to the fact that
$M_k$ are not nonnegative definite operators.
\end{remark}

\subsection{Estimation of $d$ and $\calM$}\label{section:2.2.2}

\subsubsection{Estimators and fitted dynamic models}\label{section:2.2.2.1}

Let $ \psi_1,\ldots, \psi_d $ be the orthonormal eigenfunctions of $K$ corresponding
to its $d$
nonzero eigenvalues. Then they
form an orthonormal basis of $\calM$; see Proposition~\ref{prop:1}(ii) above.
Hence, it holds that
\[
X_t(u) - \mu(u) = \sum_{j=1}^d \xi_{tj} \varphi_j(u)
= \sum_{j=1}^d \eta_{tj} \psi_j(u),
\]
where $\eta_{tj} = \int_\calI\{ X_t(u) - \mu(u) \}\psi_j(u) \,du $.
Therefore, the serial dependence of $X_t(\cdot)$ [and also that of
$Y_t(\cdot)$] can be
represented by that of
the $d$-vector process $\bfeta_t \equiv(\eta_{t1},\ldots, \eta_{td})'$.
Since $(\xi_{tj}, \varphi_j)$ cannot be estimated directly from $Y_t$
(see Section~\ref{section:2.2.1}
above), we estimate $(\eta_{tj}, \psi_j)$ instead.

As we have stated above, $M_k$ for $k\ne0$ may be directly
estimated from the observed curves $Y_t$; see (\ref{kernel}). Hence,
a natural estimator for $K$ may be defined as
%
\begin{eqnarray} \label{a35}
\widehat K(u,v) &=& \sum_{k=1}^{p} \int_\calI\widehat M_k (u, z) \widehat M_k(v, z) \,dz
\nonumber\\[-2pt]
&=& {1 \over(n-p)^2} \sum_{t,s=1}^{n-p} \sum_{k=1}^{p} \{
Y_{t}(u) -
\bar Y(u) \}\{ Y_{s}(v) - \bar Y(v) \} \\[-2pt]
&&\hspace*{77.23pt}{}\times\inner{ Y_{t+k} -\bar Y}{
Y_{s+k} - \bar Y},\nonumber
\end{eqnarray}
see (\ref{a31}), (\ref{a34}), (\ref{kernel}) and (\ref{a36}).\vadjust{\goodbreak}

By Proposition~\ref{prop:1}, we define $\widehat d$ to be the number of
nonzero eigenvalues of~$\widehat K$ (see Section~\ref{section:2.2.2.3}
below) and $\widehat\calM$ to be the linear space spanned by the $\widehat d$
corresponding orthonormal eigenfunctions $\widehat\psi_1(\cdot),
\ldots, \widehat\psi_{\widehat d}(\cdot)$. This leads to the fitting
%
\begin{equation} \label{a38}
\widehat Y_t (u) = \bar Y(u) + \sum_{j=1}^{\widehat d} \widehat\eta_{tj} \widehat
\psi_j(u),\qquad u \in\calI,
\end{equation}
where
%
\begin{equation}\label{a10:2}
\widehat\eta_{tj} = \int_\calI\{ Y_t(u) - \bar Y(u) \} \widehat\psi_j(u) \,d
u, \qquad j=1, \ldots, \widehat d.
\end{equation}
Although\vspace*{2pt} $\widehat\psi_j$ are not the estimators for
the eigenfunctions $\varphi_j$ of $M_0$ defined in (\ref{b1}),
$\widehat\calM= \operatorname{span}\{ \widehat
\psi_1(\cdot), \ldots, \widehat\psi_{d}(\cdot) \}$ is a consistent
estimator of $\calM= \operatorname{span}\{ \varphi_1(\cdot),\ldots,
\varphi_d (\cdot) \}$ (Theorem~\ref{theorem2} in Section~\ref{section:2.3} below).

In order to model the dynamic behavior of $Y_t(\cdot)$, we only need to
model the $\widehat d$-dimensional vector process $\widehat\bfeta_t \equiv(\widehat
\eta_{t1},\ldots,
\widehat\eta_{t\widehat d})'$; see (\ref{a38}) above.
This may be done using VARMA or any other multivariate time series models.
See also \citet{tiao1989} for applying linear transformations in order
to obtain a
more parsimonious model for $\widehat{\bfeta}_t$.

The integer $p$ used in (\ref{kernel}) may be selected in the same
spirit as the maximum lag used in, for example, the Ljung--Box--Pierce
portmanteau test for white noise. In practice, we often choose $p$ to be
a small positive integer. Note that $k_0$ fulfilling the
condition of Proposition~\ref{prop:1} is often small since serial
dependence decays as the lag increases
for most practical data.

\subsubsection{Eigenanalysis}\label{section:2.2.2.2}

To perform an eigenanalysis in a Hilbert space is not a trivial
matter. A popular pragmatic approach is to use an approximation via
discretization, that is, to evaluate the observed curves at a fine grid
and to replace the observed curves by the resulting vectors. This
is an approximate method;
effectively transform the problem to an eigenanalysis for a finite
matrix. See, for example, Section 8.4 of \citet{ramsay2005}. Below we also
transform the problem into an eigenanalysis of a finite matrix but
not via any approximations. Instead we make use of the well-known
duality property that $\bA\bB'$ and $\bB' \bA$ share the same
nonzero eigenvalues for any matrices $\bA$ and $\bB$ of the same
sizes. Furthermore, if $\bgamma$ is an eigenvector of $\bB'\bA$,
$\bA\bgamma$ is an eigenvector of $\bA\bB' $ with the same
eigenvalue. In fact, this duality also holds for operators in a
Hilbert space. This scheme was adopted in \citet{kneip2001} and
\citet{benko2009}.

We present a heuristic argument first.
To view the operator $\widehat K(\cdot, \cdot)$ defined
in (\ref{a35}) in the form of $\bA\bB'$, let us denote the curve
$Y_t (\cdot) -\bar Y(\cdot)$ as an $\infty\times1$ vector $\bY_t$
with $\bY_t' \bY_s = \inner{Y_t - \bar Y}{ Y_s - \bar Y}$; see
(\ref{a36}). Put\vadjust{\goodbreak} $\calY_k = (\bY_{1+k},\ldots, \bY_{n-p+k})$. Then
$\widehat K(\cdot, \cdot)$ may be represented as an $\infty\times
\infty$ matrix
\[
\widehat\bK= {1 \over(n-p)^2 } \calY_0 \sum_{k=1}^{p} \calY_k'
\calY_k \calY_0'.
\]
Applying the duality
with
$\bA= \calY_0$ and $\bB^\prime= \sum_{1 \leq k \leq p}
\calY_k' \calY_k \calY_0'$, $\widehat
\bK$ shares the same nonzero eigenvalues with the $(n-p) \times
(n-p)$ matrix
%
\begin{equation}\label{K}
\bK^* = {1 \over(n-p)^2 } \sum_{k=1}^{p} \calY_k' \calY_k
\calY_0' \calY_0,
\end{equation}
where the $(t,s)$th element of $\calY_k' \calY_k$ is
$\bY_{t+k}^\prime\bY_{s+k} = \inner{Y_{t+k} - \bar Y}{ Y_{s+k} -
\bar Y}$ and $k =0, 1,\ldots, p$. Furthermore, let $\bgamma_j=
(\gamma_{1j},\ldots, \gamma_{n-p,j})'$, $j=1,\ldots, \widehat d$, be
the eigenvectors of $\bK^*$ corresponding to the $\widehat d$ largest
eigenvalues. Then
%
\begin{equation} \label{a39}
\sum_{t=1}^{n-p} \gamma_{tj}
\{Y_t(\cdot) - \bar Y(\cdot) \},\qquad j=1,\ldots, \widehat d,
\end{equation}
are the $\widehat d$ eigenfunctions of $\widehat K(\cdot, \cdot)$. Note that
the functions in (\ref{a39}) may not be
orthogonal with each other. Thus, the orthonormal eigenfunctions $\widehat
\psi_1(\cdot),\ldots, \widehat\psi_{\widehat d}(\cdot)$ used in ({\ref{a38})
may be obtained by applying a Gram--Schmidt algorithm to the functions
given in (\ref{a39}).

The heuristic argument presented above is justified by result below.
The formal proof is relegated to Appendix~\ref{appB}.
\begin{proposition}\label{proposition2}
The operator $\widehat K(\cdot,\cdot)$ shares the same nonzero
eigenvalues with matrix $\bK^*$ defined in (\ref{K}) with the
corresponding eigenfunctions given in (\ref{a39}).
\end{proposition}
\begin{remark}\label{remark2}
The truncation of the sums in
(\ref{kernel}) at $(n-p)$ for different $k$ is necessary to ensure
the applicability of the above duality operation. If we truncated
the sum for $\widehat M_k$ at $(n-k)$ instead, $\calY_k ^\prime\calY_k $
would be of different sizes for different~$k$, and $\bK^*$ in
(\ref{K}) would not be well defined.
\end{remark}

\subsubsection{Determination of $d$ via statistical tests}
\label{section:2.2.2.3}

Although the number of nonzero eigenvalues of operator $K(\cdot,
\cdot)$ defined in~(\ref{a31}) is $d$ [Proposition~\ref{prop:1}(ii)], the number of nonzero eigenvalues of its
estimator $\widehat K(\cdot, \cdot)$ defined in~(\ref{a35}) may be
much greater than $d$ due to random fluctuation in the sample. One
empirical approach is to take $\widehat d$ to be the number of
``large''
eigenvalues of $\widehat K$ in the sense that the $(\widehat d+1)$th largest
eigenvalue drops significantly; see also Theorem~\ref{theorem3} in
Section~\ref{section:2.3} and
Figure~\ref{figure:1} in Section~\ref{section:2.4.1}.
Hyndman and Ullah (\citeyear{hyndman2007}) proposed to choose $d$ by minimizing forecasting
errors.
Below, we present a bootstrap test to determine the value of $d$.\vadjust{\goodbreak}

Let $\theta_1 \ge\theta_2 \ge\cdots\ge0$ be the eigenvalues
of $K$.
If the true dimensionality is $d=d_0$, we expect to reject the null hypothesis
$\theta_{d_0} =0$, and not to reject the hypothesis $\theta_{d_0+1}=0$.
Suppose we are interested in testing the null hypothesis
%
\begin{equation} \label{h0}
H_0\dvtx \theta_{d_0+1} =0,
\end{equation}
where $d_0$ is a known integer, obtained, for example, by visual
observation of the estimated eigenvalues $\widehat\theta_1 \ge\widehat
\theta_2 \ge\cdots\ge0$ of $\widehat K$. Hence, we reject $H_0$ if
$\htheta_{d_0 +1} > l_{\alpha}$, where $l_{\alpha}$ is the critical
value at the $\alpha\in(0,1)$ significance level. To evaluate the
critical value $l_\alpha$, we propose the following bootstrap
procedure.

\begin{enumerate}
\item Let $\widehat Y_t(\cdot)$ be defined as in (\ref{a38}) with
$\widehat d = d_0$. Let $\widehat\ve_t(\cdot) = Y_t(\cdot) - \widehat Y_t(\cdot)$.
\item Generate a bootstrap sample from the model
\[
Y_t^*(\cdot) = \widehat Y_t(\cdot) + \ve_t^*(\cdot),
\]
where $ \ve_t^*$ are drawn independently (with replacement) from $\{
\widehat\ve_1,\ldots, \widehat\ve_n\}$.
\item
Form an operator $K^*$ in the same manner as $\widehat K$ with $\{ Y_t\}$ replaced
by $\{ Y_t^* \}$, compute the $(d_0+1)$th largest
eigenvalue $\theta^*_{d_0+1}$ of $K^*$.
\end{enumerate}
Then the conditional distribution of $\theta^*_{d_0+1}$, given the
observations $\{ Y_1,\ldots, Y_n \}$, is taken as the distribution
of $\htheta_{d_0 +1}$ under $H_0$. In practical implementation, we
repeat Steps 2 and 3 above $B$ times for some large integer $B$, and
we reject $H_0$ if the event that $\theta^*_{d_0+1}> \htheta_{d_0
+1}$ occurs not more than $[\alpha B]$ times. The simulation results
reported in Section~\ref{section:2.4.1} below indicate that the
above bootstrap method works well.
\begin{remark}\label{remark3}
The serial dependence in $X_t$ could provide an alternative method for testing
hypothesis (\ref{h0}).
Under model (\ref{a5}), the projected series
of the curves $Y_t(\cdot)$ on any direction perpendicular to $\calM$
is white noise.
Put $U_t = \inner{ Y_t}{\widehat\psi_{d_0+1}}, t=1,\ldots, n$.
Then $U_t$ would behave like a (scalar) white noise under $H_0$.
However, for example, the Ljung--Box--Pierce portmanteau test for white
noise coupled with the standard $\chi^2$-approximation does not
work well in this context.
This is due to the fact that the $(d+1)$th largest eigenvalue $\widehat K$
is effectively the
extreme value of the estimates for all the zero-eigenvalues of $K$.
Therefore, $\widehat\psi_{d_0+1}$ is not an estimate for a fixed direction,
which makes the
$\chi^2$-approximation for the Ljung--Box--Pierce statistic
mathematically invalid.
Indeed some simulation results, not reported here, indicate that
the $\chi^2$-approximation tends to
underestimate the critical values for the Ljung--Box--Pierce test in
this particular context.
\end{remark}


\section{Theoretical properties}\label{section:2.3}


Before presenting the asymptotic results, we first solidify some
notation. Denote by $(\theta_j, \psi_j)$ and $(\widehat\theta_j,
\widehat\psi_j)$ the (eigenvalue, eigenfunction) pairs of $K$ and $\widehat K$,
respectively [see (\ref{a31}) and (\ref{a35})]. We always arrange
the eigenvalues\vspace*{1pt} in descending order, that is,
$\theta_j > \theta_{j+1}$. As the eigenfunctions of $K$ and $\widehat K$
are unique only up to sign changes, in the sequel, it will go without
saying that the right versions are used. Furthermore, recall that
$\theta_j = 0$ for all $j \geq d+1$. Thus, the eigenfunctions
$\psi_{j}$ are not identified for $j \geq d+1$. We take this last
point into consideration in our theory.
We always assume that the dimension $d \ge1$ is a fixed finite
integer, and $p\ge1$ is also a fixed finite integer.

For simplicity in the proofs, we suppose that $E\{Y_t(\cdot)\} = \mu
(\cdot)$ is known and thus set $\bar{Y}(\cdot) = \mu(\cdot)$.
Straightforward adjustments to our arguments can be made when this is
not the case. We denote by $\hilbert{L}$ the Hilbert--Schmidt norm for
any operator~$L$; see Appendix~\ref{appA}. Our asymptotic results are based on
the following regularity conditions:

\begin{longlist}[(C3)]

\item[C1.] $\{ Y_t(\cdot) \}$ is strictly stationary and $\psi$-mixing with
the mixing coefficient defined as
\[
\psi(l) = \sup_{A \in\calF^0_\infty, B \in\calF^\infty_l,
P(A)P(B) > 0} |1 - P(B|A)/P(B)|,
\]
where $\calF_i^j = \sigma\{ Y_i(\cdot),\ldots, Y_j(\cdot) \}$ for any
$j\geq i$.
In addition, it holds that $\sum_{l=1}^\infty l\times\psi^{1/2}(l) <
\infty$.\vspace*{1pt}

\item[C2.] $E\{\int_\calI Y_t(u)^2 \,du\}^2 < \infty$.
\item[C3.] $\theta_1 > \cdots> \theta_d > 0 = \theta_{d+1} = \cdots,$
that is, all the nonzero
eigenvalues of $K$ are different.
\item[C4.] $\cov\{X_s(u), \ve_t(v) \} = 0$ for all $s,t$ and $u,v \in
\calI$.
\end{longlist}
\begin{theorem}\label{theorem1}
Let
conditions \textup{C1--C4} hold. Then as $n\to\infty$, the following assertions
hold:

\begin{longlist}
\item $\hilbert{\widehat K - K} = O_p(n^{-1/2})$.

\item For $j=1,\ldots,d$, $| \widehat\theta_j - \theta_j | =
O_p(n^{-1/2})$ and
\[
\biggl(\int_\calI \{ \widehat\psi_j(u) - \psi_j(u) \}^2 \,du
\biggr)^{1/2} = O_P(n^{-1/2}).
\]

\item For $j \ge d+1$, $\widehat\theta_j= O_p(n^{-1})$.

\item Let $\{\psi_j \dvtx j \geq d+1\}$ be a complete orthonormal basis of
$\calM^\bot$,
and put
\[
f_j(\cdot) = \sum_{i=d+1}^\infty\inner{\psi_i}{\widehat\psi_{j}}
\psi_i(\cdot).
\]
Then for any $j\ge d+1$,
\[
\Biggl( \int_\calI \Biggl\{\sum_{i=1}^d \inner{\psi_i}{\widehat\psi_{j}}
\psi_i(u) \Biggr\}^2 \,du \Biggr)^{1/2} = \biggl( \int_\calI \{ \widehat
\psi_{j}(u) - f_j(u) \}^2 \,du \biggr)^{1/2} =
O_p(n^{-1/2}).
\]
\end{longlist}
\end{theorem}
\begin{remark}\label{remark4}
(a) In the above theorem,
assertions (i) and (ii) are standard.
(In fact, those results still hold for $d=\infty$.)\vadjust{\goodbreak}

(b) Assertion (iv)
implies that the estimated eigenfunctions $\widehat\psi_{d+j}$, $j \ge1$,
are asymptotically in the orthogonal complement of the
dynamic space $\calM$.

(c) The fast convergence rate $n$ in assertion (iii) deserves some
further explanation.
To this end, we consider a simple analogue:
let $A_1,\ldots,A_n$ be a sample of stationary random variables, and
we are interested in estimating $\mu^2 = (EA_t)^2$ for which we use the
estimator $\bar{A}^2 = (n^{-1} \sum_{t=1}^n A_t)^2 = n^{-2}\sum_{s,t=1}^n
A_s A_t$. Then under appropriate regularity conditions, it holds that
%
\begin{equation}\label{example}\qquad
| \bar{A}^2 - \mu^2 | \leq|\mu| |\bar{A} - \mu| + |\bar{A}^2 - \bar{A}
\mu| = |\mu| \cdot O_p(n^{-1/2}) + O_p(n^{-1})
\end{equation}
as $ |\bar{A} - \mu| = O_p(n^{-1/2})$
and $ |\bar{A}^2 - \bar{A} \mu| = O_p(n^{-1})$. The latter
follows\break
from~a~simple~$U$-statistic argument; see \citet{lee1990}. It is easy to see\break from
(\ref{example}) that
$| \bar{A}^2 - \mu^2 | = O_p(n^{-1/2})$ if $\mu\neq0$, and
$| \bar{A}^2 - \mu^2 | = O_p(n^{-1})$ if\break $\mu= 0$.
In our context, the operator $\widehat K = \sum_{k=1}^p \int_\calI M_k(u, r)
M_k(v,r) = \break(n-p)^{-2}
\sum_{k=1}^p \sum_{s,t=1}^{n-p} Z_{ik} Z_{jk}^*(u,v)$, where
$Z_{tk}(u,v) = \{Y_t(u) - \mu(u)\} \{Y_{t+k}(v) - \mu(v)\}$ and $Z_{ik}
Z_{jk}^*(u,v) = \int_\calI Z_{ik}(u, r) Z_{jk}(v,r) \,dr$, is
similar to $\bar{A}^2$, and hence the convergence properties stated in
Theorem~\ref{theorem1}(iii) [and also (ii)].
The fast convergence rate, which is termed as
``super-consistent'' in econometric literature,
is illustrated via simulation
in Section~\ref{section:2.4.1} below; see Figures~\ref{figure:4}--\ref{figure:7}.
It makes the identification of zero-eigenvalues
easier; see Figure~\ref{figure:1}.

\begin{figure}

\includegraphics{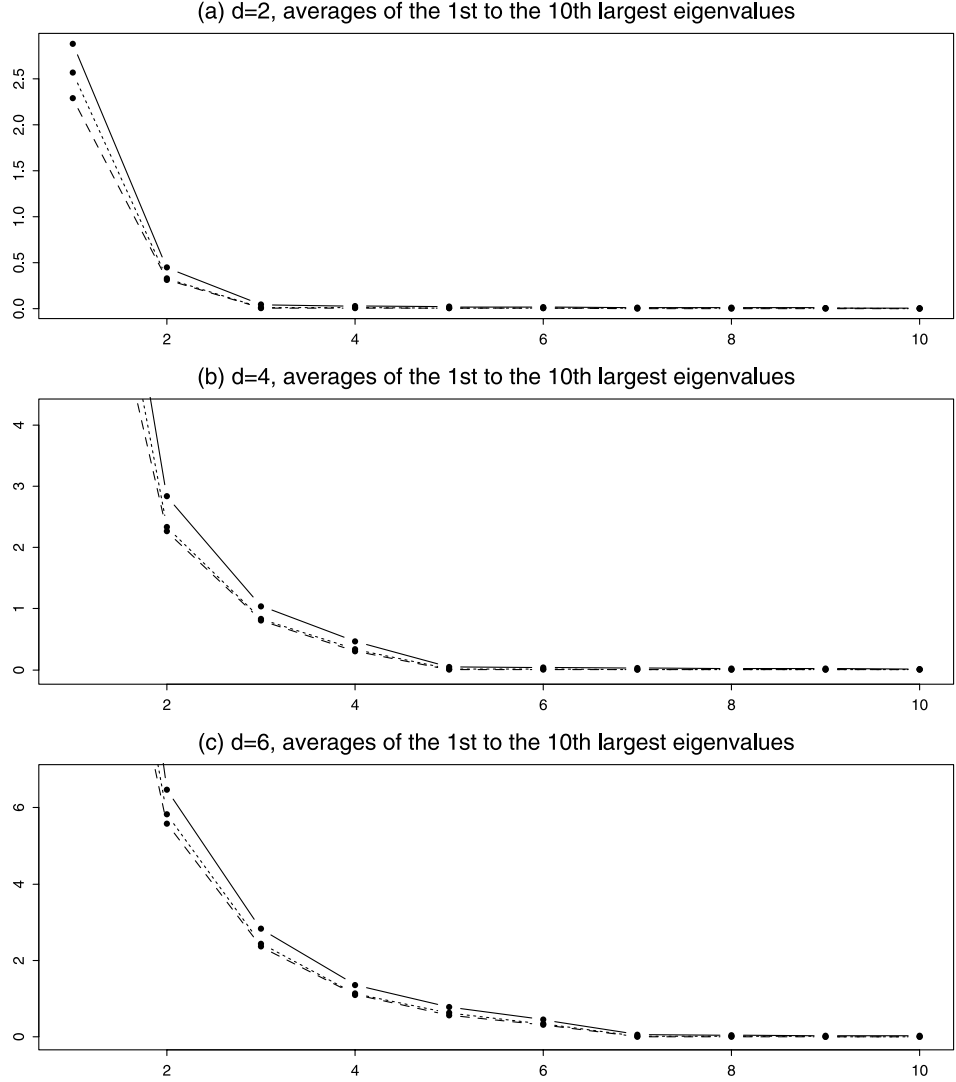}%
\vspace*{-3pt}
\caption{The average estimated eigenvalues over the 200 replications with
sample sizes $n=100$ (solid lines), $300$ (dotted lines) and $600$
(dashed lines).}\label{figure:1}\vspace*{-3pt}
\end{figure}

With $d$ known, let $\widetilde
\calM= \operatorname{span} \{ \widehat\psi_1(\cdot),\ldots, \widehat\psi_{
d}(\cdot)\}$, where $\widehat\psi_1(\cdot),\ldots, \widehat\psi_{
d}(\cdot)$ are the eigenfunctions of $\widehat K$ corresponding to the
$d$ largest eigenvalues. In order
to measure the discrepancy between
$\calM$ and $\widetilde\calM$, we introduce the following metric. Let $\calN_1$
and $\calN_2$ be any
two $d$-dimensional subspaces of $\calL_2(\calI)$. Let $\{
\zeta_{i1}(\cdot),\ldots, \zeta_{id}(\cdot) \}$ be an orthonormal
basis of $\calN_i$, $i=1, 2$. Then the projection of $\zeta_{1k}$
onto $\calN_2$ may be expressed as
\[
\sum_{j=1}^d \inner{\zeta_{2j}}{\zeta_{1k}} \zeta_{2j}(u).
\]
Its squared norm is $ \sum_{j=1}^d (
\inner{\zeta_{2j}}{\zeta_{1k}})^2 \le1. $ The discrepancy measure
is defined as
%
\begin{equation} \label{metric}
D(\calN_1, \calN_2) = \sqrt{ 1 - {1 \over d} \sum_{j,k=1}^d (
\inner{\zeta_{2j}}{\zeta_{1k}})^2}.
\end{equation}
It is clear that this is a symmetric measure between 0 and 1. It is
independent of the choice of the orthonormal bases used in the
definition, and it equals 0 if and only if $\calN_1=\calN_2$. Let
$\calZ$ be the set consisting of all the $d$-dimensional subspaces
in $\calL_2(\calI)$. Then $(\calZ, D)$ forms a metric space in the
sense that $D$ is a well-defined distance measure on $\calZ$ (see
Lemma~\ref{lemma4} in Appendix~\ref{appB} below).
\end{remark}
\begin{theorem}\label{theorem2}
Let the conditions of Theorem~\ref{theorem1}   hold. Suppose that
$d$ is known.
Then as $n \to\infty$, it holds that $D(\widetilde{\calM},\calM) =
O_p(n^{-1/2})$.\vadjust{\goodbreak}
\end{theorem}
\begin{remark}\label{remark5}
Our estimation of $\calM$ is asymptotically
adaptive to $d$.
To this end, let $\widehat d$ be a consistent estimator of $d$ in the sense
that $P(\widehat d = d) \to1$, and $\widehat\calM=
\operatorname{span} \{ \widehat\psi_1, \ldots, \widehat\psi_{\widehat d}\}$ be the estimator
of $\calM$ with $d$ estimated by $\widehat d$.
Since $\widehat d$ may differ
from $d$, we use the modified metric $\widetilde D$, defined in (\ref
{metric4}) below,
to measure the difference between $\widehat\calM$ and $\calM$. Then it
holds for any
constant $C>0$ that
\begin{eqnarray*}
&& P \{ n^{1/2} | \widetilde D(\widehat\calM, \calM) - D(\widetilde\calM, \calM) | > C
\} \\
&&\qquad\leq P \{ n^{1/2} |\widetilde D(\widehat\calM, \calM) - D(\widetilde\calM, \calM) | > C
| \widehat d = d \} P(\widehat d = d) + P(\widehat d \neq d)\\
&&\qquad\leq P \{ n^{1/2} |\widetilde D(\widehat\calM, \calM) - D(\widetilde\calM, \calM) | > C
| \widehat d = d \} + o(1).
\end{eqnarray*}
Note that when $\widehat d = d$,
$\widehat\calM= \widetilde\calM$ and thus $\widetilde D(\widehat\calM, \calM) = D(\widehat\calM,
\calM)$.
Hence the conditional probability on the RHS of the above expression is 0.
This together with Theorem~\ref{theorem2} yield $\widetilde D(\widehat\calM, \calM)
= O_p(n^{-1/2})$.

One such consistent estimator of $d$ may be defined as $\widehat{d} = \#
\{j \dvtx \widehat\theta_j \geq\epsilon\}$, where $\epsilon= \epsilon(n) >0$
satisfies the
conditions in Theorem~\ref{theorem3} below.
\end{remark}
\begin{theorem}\label{theorem3}
Let the conditions of Theorem~\ref{theorem1}   hold. Let $\epsilon
\to0$ and
$\epsilon^2 n \to\infty$ and as $n \to\infty$.
Then $P( \widehat{d} \ne d) \to0$.
\end{theorem}


\section{Numerical properties}\label{section:2.4}
\subsection{Simulations}\label{section:2.4.1}

We illustrate the proposed method first using the simulated data from
model ({\ref{a1}) with
\[
X_t(u) = \sum_{i=1}^d \xi_{ti}\varphi_i(u),\qquad \ve_t(u) =
\sum_{j=1}^{10} {Z_{tj} \over2^{j-1}} \zeta_j(u),\qquad u\in
[0,1],
\]
where $\{ \xi_{ti}, t\ge1 \}$ is a linear AR(1) process with the
coefficient $(-1)^i(0.9-0.5i/d)$, the innovations $Z_{tj}$ are
independent $N(0, 1)$ variables and
\[
\varphi_i(u) = \sqrt{2} \cos(\pi i u),\qquad
\zeta_j(u) = \sqrt{2} \sin(\pi ju).
\]
We set sample size $n=100, 300 $ or 600, and the dimension
parameter $d =2, 4 $ or~6. For each setting, we repeat the
simulation 200 times. We use $p=5$ in defining the operator $\widehat K$
in (\ref{a35}). For each of the 200 samples, we replicate
the bootstrap sampling 200 times.

The average of the ordered eigenvalues of $\widehat K$ obtained from the
200 replications are plotted in Figure~\ref{figure:1}. For a good visual
illustration, we only plot the ten largest eigenvalues. It is clear
that drop from the $d$th largest eigenvalue to the $(d+1)$st is very
pronounced.
Furthermore, the estimates for zero-eigenvalues with different sample
size are much closer than those for nonzero eigenvalues. This evidence
is in line with the different convergence rates presented in
Theorem~\ref{theorem1}(ii) and (iii).
We apply the bootstrap method to test the hypothesis
that the $d$th or the $(d+1)$st largest eigenvalue of $K$
($\theta_d$ and $\theta_{d+1}$, resp.) are 0. The results are
summarized in Figure~\ref{figure:2}. The
bootstrap test cannot reject the true null hypothesis
$\theta_{d+1}=0$. The false null hypothesis $\theta_d =0$ is
routinely rejected when $n=600$ or 300; see
Figure~\ref{figure:2}(a). However, the test does not
work when the sample size is as small as~100.

\begin{figure}

\includegraphics{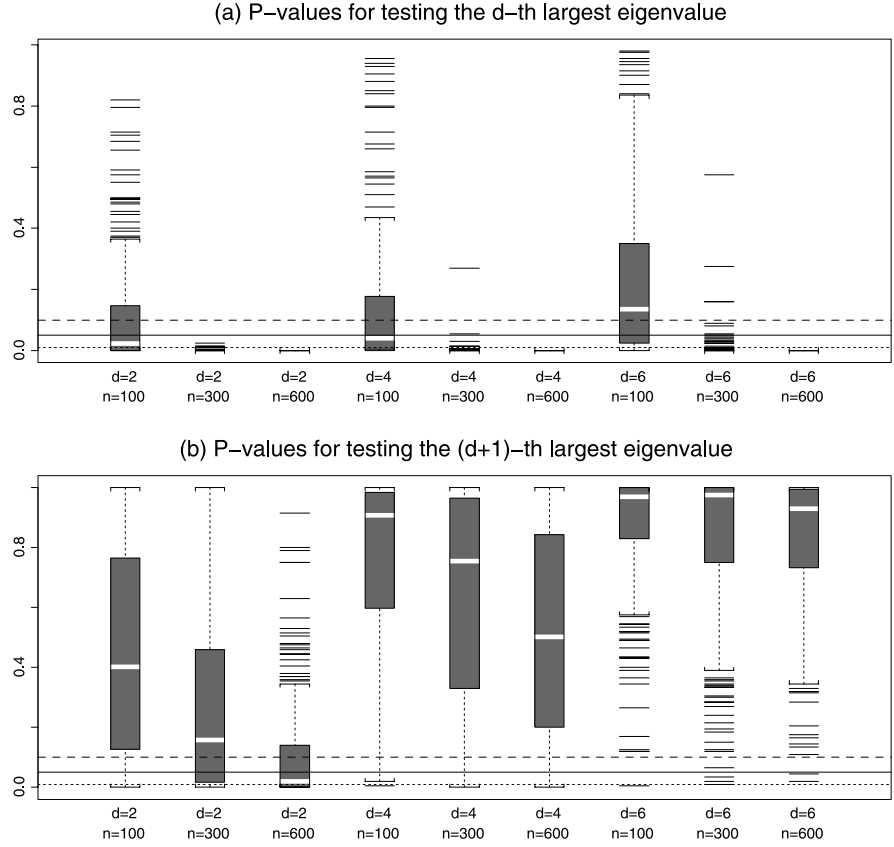}

\caption{The boxplots of the $P$-values for the bootstrap tests of the
hypothesis that \textup{(a)} the $d$th largest eigenvalue of $K$ is 0, and
\textup{(b)} the $(d+1)$th largest eigenvalue of $K$ is 0. The horizontal lines
mark the 1\% (dotted line), 5\% (solid lines) and 10\% (dashed lines)
significance levels, respectively.}\label{figure:2}
\end{figure}

To measure the accuracy of the estimation for the factor loading
space $\calM$, we need to modify the metric $D$ defined in
(\ref{metric}) first, as $\widehat d $ may be different from $d$.
Let $\calN_1, \calN_2$ be two subspaces in $\calL_2(\calI)$
with\vadjust{\goodbreak}
dimension $d_1$ and $d_2$, respectively. Let $\{ \zeta_{i1},\ldots,
\zeta_{id_i} \}$ be an orthonormal basis of $\calN_i$, $i=1,2$. The
discrepancy measure between the two subspaces is defined as
%
\begin{equation} \label{metric4}
\widetilde D(\calN_1, \calN_2) =
\sqrt{ 1 - {1 \over\max( d_1, d_2) } \sum_{k=1}^{d_1} \sum_{j=1}^{d_2}
(\inner{\zeta_{2j}}{\zeta_{1k}})^2 }.
\end{equation}
It can be shown that $\widetilde D(\calN_1, \calN_2) \in[0, 1]$.
It equals 0 if and only if $\calN_1 = \calN_2$, and~1 if and only if
$\calN_1\,\bot\,\calN_2$.
Obviously, $\widetilde D(\calN_1, \calN_2) = D(\calN_1, \calN_2)$
when $d_1 = d_2 = d$.
We computed $\widetilde D(\widehat\calM, \calM)$ in the 200 replications
for each setting.
Figure~\ref{figure:3} presents the boxplots of those
$\widetilde D$-values. It is noticeable that the $\widetilde D$
measure decreases as the sample size $n$ increases. It is
interesting to note too that the accuracy of the estimation is
independent of the dimension $d$.

\begin{figure}

\includegraphics{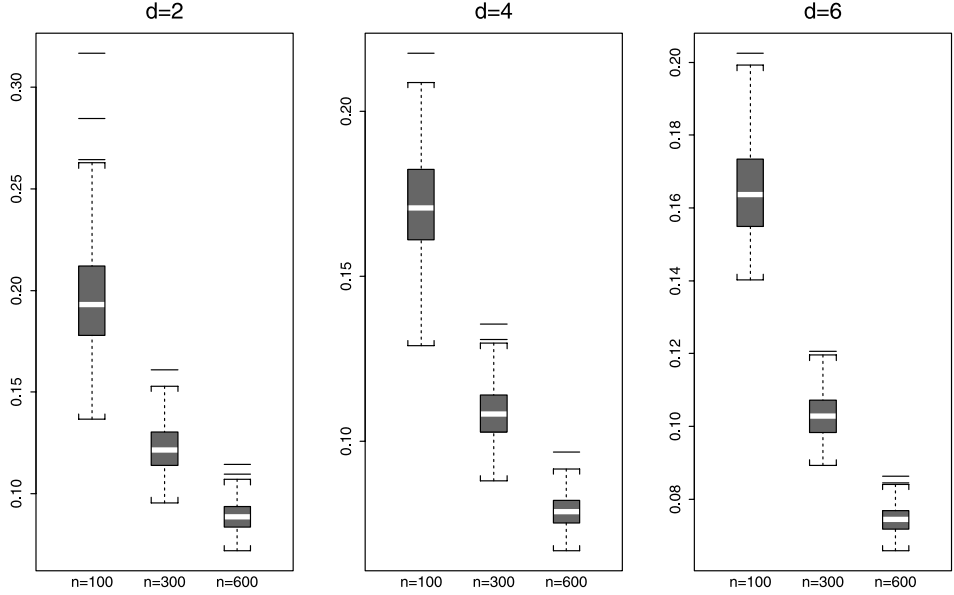}

\caption{The boxplots for the estimated error $\widetilde D$ defined in
(\protect\ref{metric4}).}\label{figure:3}\vspace*{-3pt}
\end{figure}

To further illustrate the different convergence rates in estimating
nonzero and zero eigenvalues, as stated in Theorem~\ref{theorem1}, we generate
10,000 samples with different sample sizes from
model (\ref{a1})
with $d=1$, $\xi_t = 0.5 \xi_{t-1} + \eta_t$, where $\eta_t \sim
N(0,1)$, $\varphi(u) = \sqrt{2} \cos(\pi u)$, and $\ve_t(\cdot)$ is
the same as above.
In defining the operator $K$, we let $p=1$.
Then the operator $K$ has only one nonzero
eigenvalue $\theta= 2$.
Figure~\ref{figure:4} depicts the standardized histograms and the
kernel density estimators of $\sqrt{n} (\widehat
\theta_1 - \theta)$, computed from the 10,000 samples.
It is evident that those distributions resemble normal distributions
when the sample size is 200 or greater.
This is in line with
Theorem~\ref{theorem1}(ii) which implies that $\sqrt{n} (\widehat
\theta_1 - \theta)$ converges to a nondegenerate distribution.

Figure~\ref{figure:5} displays the distribution of $\sqrt{n}
\widehat\theta_2$, noting $ \theta_2=0$. It is clear that $\sqrt{n} \widehat
\theta_2$ converges to zero as $n$ increases,
indicating the fact that
the normalized factor $\sqrt{n}$ is too small to stabilize the distribution.
In contrast, Figure~\ref{figure:6} exhibits that
the distribution of $n \widehat\theta_2$
stabilizes from the sample size as small as $n=50$; see
Theorem~\ref{theorem1}(iii).
In fact, the profile of the distribution with
$n=10$ looks almost the same as that with $n=2000$.

\begin{figure}

\includegraphics{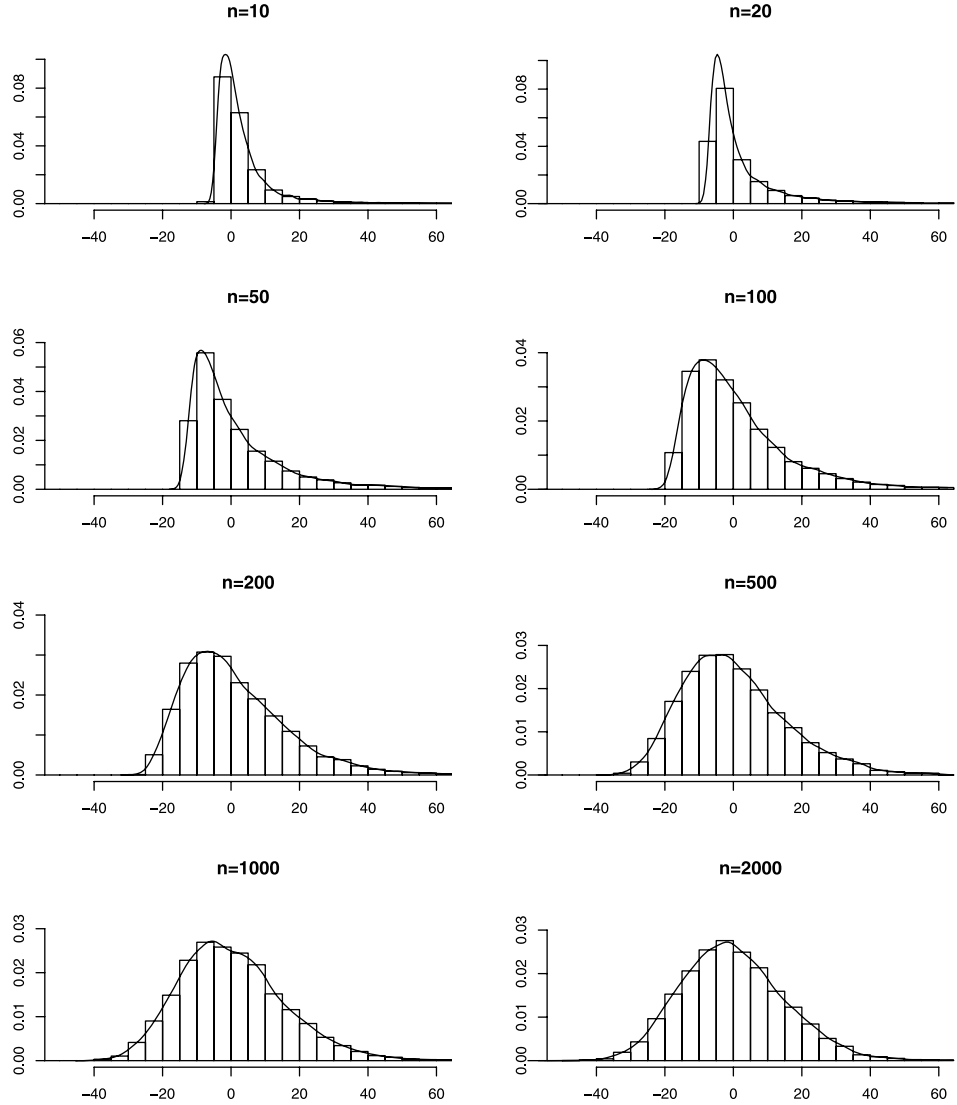}

\caption{Standardized histograms overlaid by kernel density
estimators of $\sqrt{n} (\widehat\theta_1 - \theta)$.}\label{figure:4}
\end{figure}

\begin{figure}

\includegraphics{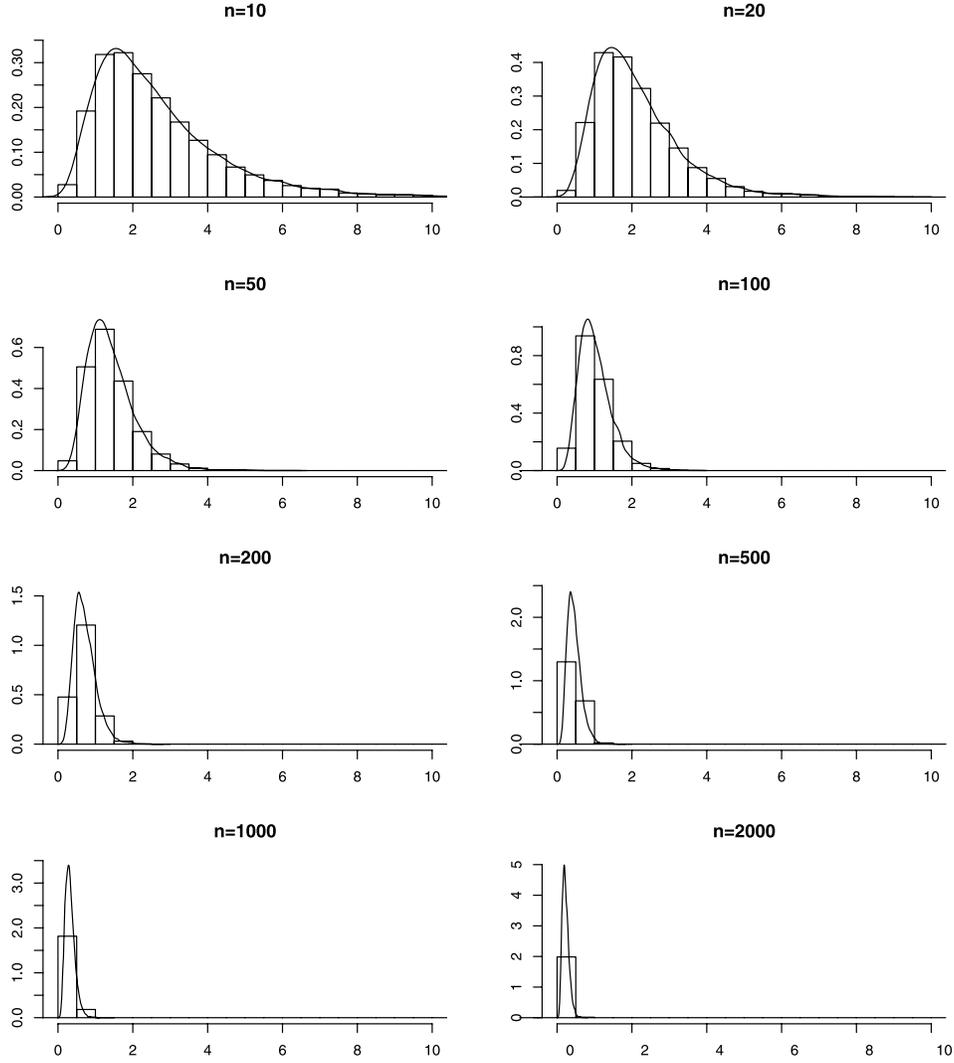}

\caption{Standardized histograms overlaid by kernel density
estimators of $\sqrt{n} \widehat\theta_2$.}\label{figure:5}\vspace*{6pt}
\end{figure}

\begin{figure}

\includegraphics{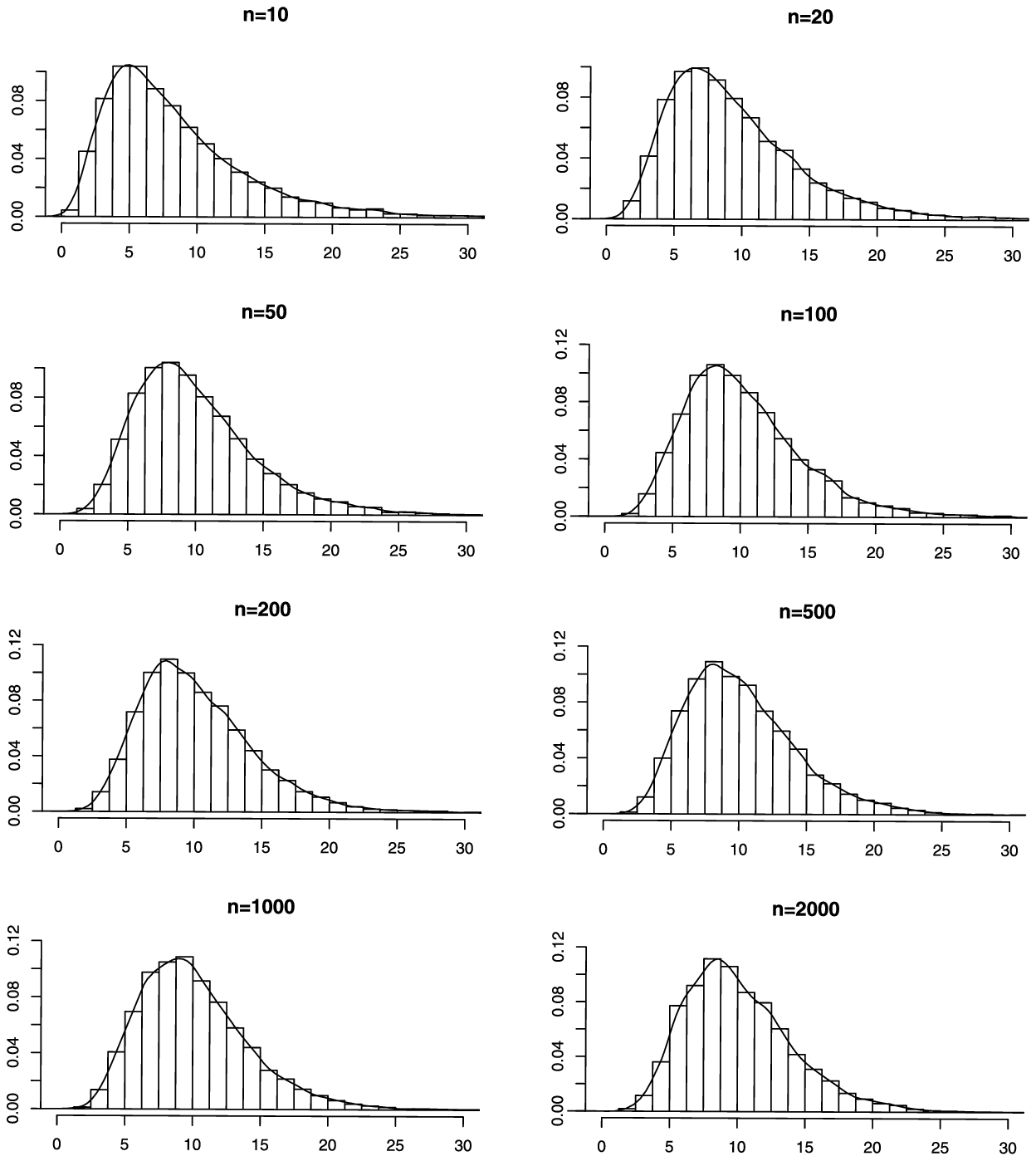}

\caption{Standardized histograms overlaid by kernel density estimators
of $n \widehat
\theta_2$.}\label{figure:6}
\end{figure}

\begin{figure}

\includegraphics{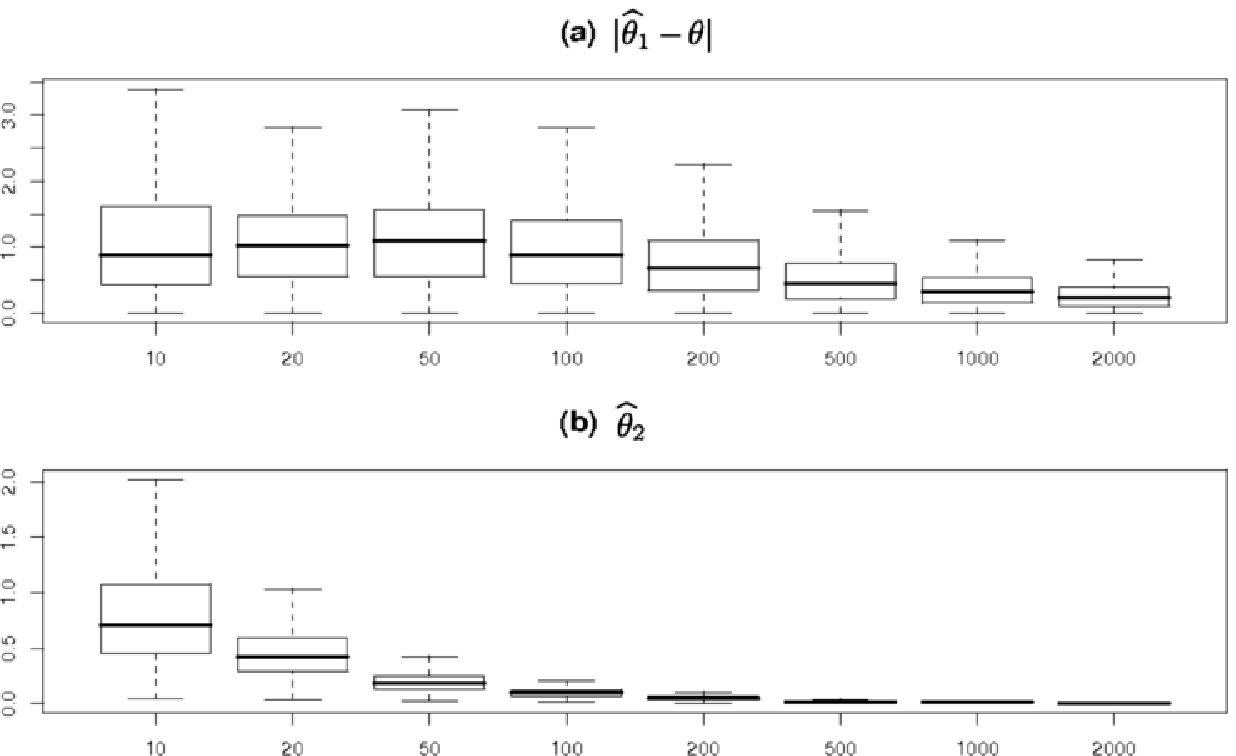}

\caption{Boxplots of estimation errors:
\textup{(a)} Errors for nonzero eigenvalue $|\widehat\theta_1 - \theta|$;
\textup{(b)}~Errors for zero-eigenvalue $\widehat\theta_2$.
To add clarity to the display, the outliers are not
plotted.}\label{figure:7}
\end{figure}

Figure~\ref{figure:7}
displays boxplots of the absolute estimation errors of the
eigenvalues. With the same sample size, the estimation errors for the
nonzero eigenvalue are considerably greater that those for the
zero eigenvalue.

\subsection{A real data example}\label{section:2.4.2}

To further illustrate the methodology developed in this paper, we set upon
the task of modeling the intraday return densities\vadjust{\goodbreak} for the IBM stock
in 2006. To this end, we have obtained the intraday prices via the
WRDS database. We only use prices between 09:30--16:00 since
the market is not particularly active outside of these times. There
are $n=251$ trading days in the sample and a total of 2,786,650
observations. The size of this dataset is 73.7~MB.

Since high frequency prices are not equally spaced in time,
we compute the returns using the prices at the so-called previous tick times
in every 5 minute intervals.
More precisely, we set the sampling times at $\tau_1 = 09$:35,
$\tau_2 = 09$:40$, \ldots, \tau_m = 16$:00 with $m=78$.
Denote by $X_{i} (t_{ij})$ the stock price on the $i$th day at the
time $t_{ij}$, $j=1,\ldots,n_i$ and $i=1,\ldots,n$.
The previous tick times on the $i$th day are defined as
\[
\tau_{il} = \max\{ t_{ij} \dvtx t_{ij} \le\tau_l, j=1,\ldots,n_i
\},\qquad
l = 1,\ldots,m.
\]
The $l$th return on the $i$th day is then defined as
${Z}_{il} = \log\{ {X}_i(\tau_{il}) / {X}_i(\tau_{i,l-1})\}$.

We then estimate the intraday
return densities using the standard kernel method
%
\begin{equation}\label{a10:1}
Y_i(u) = (m h_i)^{-1} \sum_{j=1}^m K \biggl( \frac{Z_{ij} - u}{h_i}
\biggr),\qquad i=1,\ldots,n,
\end{equation}
where $K(u) = (\sqrt{2\pi})^{-1} \exp(-u^2/2)$ is a Gaussian kernel
and $h_i$ is a bandwidth. We set $\calI= [-0.002, 0.002]$ as the support
for $Y_i(\cdot)$.
Let $\widehat\sigma_i$ be the sample standard deviation of $\{Z_{ij}, j=1,\ldots, m\}$ and $\widehat h_i = 1.06 \widehat\sigma_i m^{-1/5}$ be Silverman's
rule of thumb bandwidth choice for day $i$. Then for each $i$, we employ
three levels of smoothness by setting $h_i$ in (\ref{a10:1}) equal to $0.5
\widehat h_i$, $\widehat h_i$ and $2 \widehat h_i$.
Figure~\ref{figure:8} displays the observed densities for the first
8 days of the sample.

\begin{figure}

\includegraphics{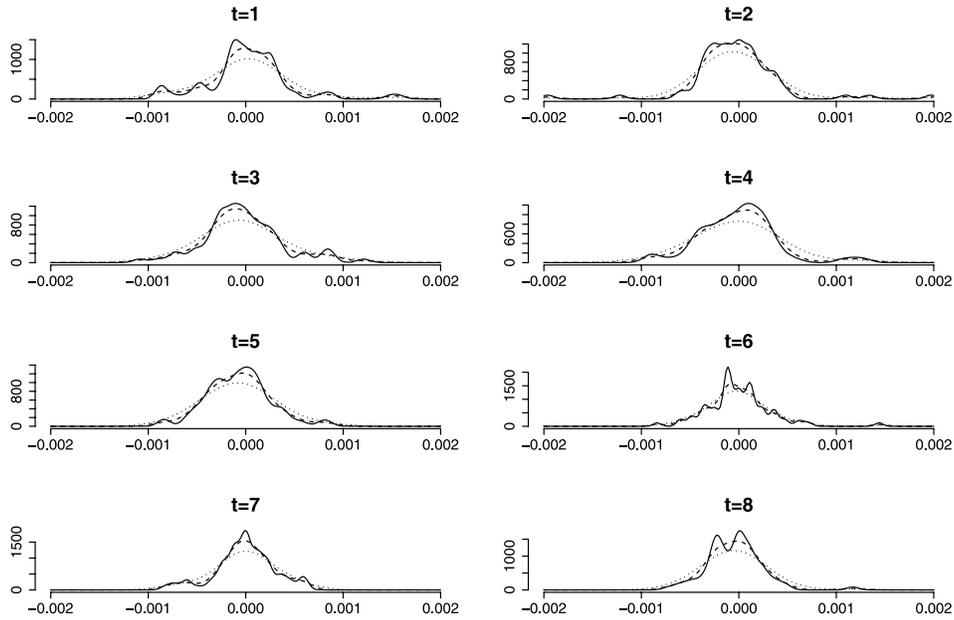}

\caption{Estimated densities, $Y_i(\cdot)$, using
bandwidths $h_i = \widehat h_i$ (solid lines), $0.5 \widehat h_i$ (dashed lines) and
$2 \widehat h_i$ (dotted lines).}\label{figure:8}
\end{figure}

\begin{figure}

\includegraphics{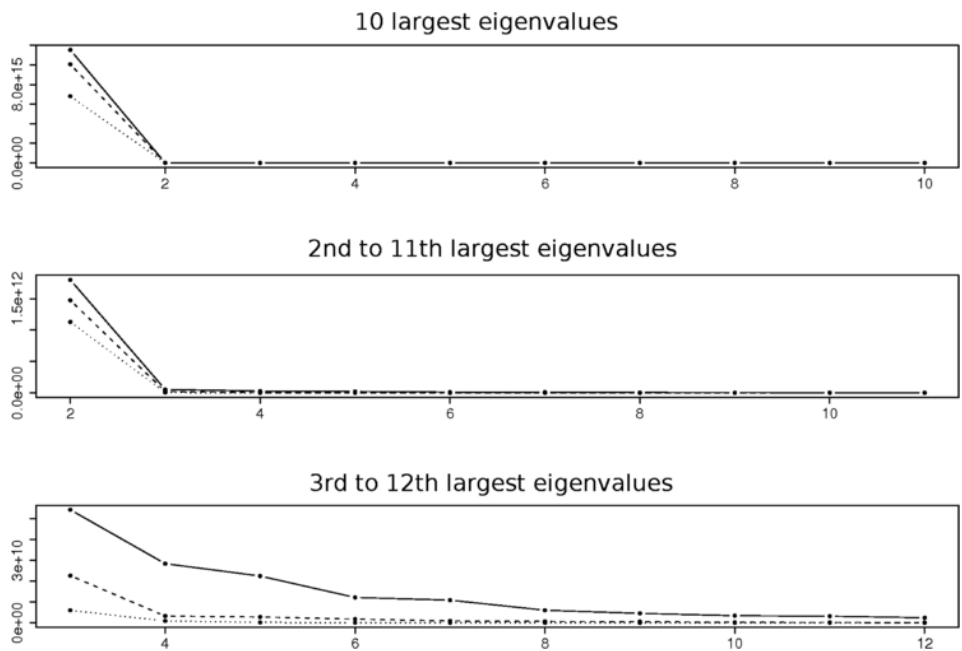}

\caption{Estimated eigenvalues $\widehat\theta_j$ using
bandwidths $h_t = 0.5 \widehat h_t$ (solid lines), $\widehat h_t$ (dashed lines) and
$2 \widehat h_t$ (dotted lines).}\label{figure:9}
\end{figure}

To identify the finite dimensionality of $Y_t(\cdot)$, we apply the methodology
developed in this paper.
We set $p=5$ in
(\ref{a31}).
Figure~\ref{figure:9} displays the estimated eigenvalues.
With all three bandwidths used, the first two eigenvalues
are much larger than the remaining ones. Furthermore, there is no
clear cut-off from the third
eigenvalue onwards. This
suggests to take $\widehat d=2$. The bootstrap tests, reported in Table \ref
{table:2},
lend further support to this assertion.
Indeed for all levels of smoothness adopted, the bootstrap test rejects
the null $H_0 \dvtx \theta_2 = 0$ but cannot reject the hypothesis $\theta_j=0$
for $j=3, 4$ or 5.
Note that it is implied by $\theta_3 = 0$ that $\theta_{3+k}=0$ for
$k\ge
1$. Indeed, we tested $\theta_{3+k}
= 0$ only for illustrative purposes.

Table~\ref{table:3} contains the $P$-values from testing the
hypothesis that the estimated loadings, $\widehat\eta_{tj}$ in
(\ref{a10:2}) are white noise using the Ljung--Box--Pierce
portmanteau test. Although we should interpret the results of this
test with caution (see Remark~\ref{remark3} in Section
\ref{section:2.2.2.3}), they provide further evidence that there is
a considerable amount of dynamic structure in the two-dimensional subspace
corresponding to the first two eigenvalues $\theta_1$ and $\theta_2$, and
there is little or none dynamic structure in the directions
corresponding to $\theta_3$ and $\theta_4$. Collating all the relevant
findings, we
comfortably set $\widehat d = 2$ in our analysis.

\begin{table}
\tablewidth=240pt
\def\arraystretch{0.9}
\caption{$P$-values from applying the bootstrap test in Section
\protect\ref{section:2.2.2.3} to the intraday
return density example}\label{table:2}
\begin{tabular*}{\tablewidth}{@{\extracolsep{\fill}}lccc@{}}
\hline
& $\bolds{h_t=0.5 \widehat h_t}$ & $\bolds{h_t = \widehat h_t}$
& $\bolds{h= 2 \widehat h_t}$ \\
\hline
$H_0\dvtx\theta_1=0$ & 0.00 & 0.00 & 0.00 \\
$H_0\dvtx\theta_2=0$ & 0.00 & 0.00 & 0.00 \\
$H_0\dvtx\theta_3=0$ & 0.35 & 0.15 & 0.18 \\
$H_0\dvtx\theta_4=0$ & 0.62 & 0.73 & 0.74 \\
$H_0\dvtx\theta_5=0$ & 0.68 & 0.91 & 0.93 \\
\hline
\end{tabular*}
\end{table}

\begin{table}[b]
\def\arraystretch{0.9}
\caption{$P$-values from testing the hypothesis $H_0 \dvtx
\widehat\eta_{tj}$ is white noise using the Ljung--Box--Pierce portmanteau
test. The test statistic is given by $Q_j = n(n+2) \sum_{k=1}^q
s_j(k)^2 / (n-k)$, where~$s_j(k)$~is the sample autocorrelation of
$\widehat\eta_{tj}$ at lag $k$. Under $H_0$, $Q_j$ has~an~asymptotic
$\chi^2_q$-distribution}\label{table:3}
\begin{tabular*}{\tablewidth}{@{\extracolsep{\fill}}lccccccccc@{}}
\hline
$\bolds{h_t}$ & \multicolumn{3}{c}{$\bolds{0.5 \widehat h_t}$} & \multicolumn{3}{c}{$\bolds{\widehat h_t}$}
& \multicolumn{3}{c@{}}{$\bolds{2 \widehat h_t}$}  \\[-4pt]
& \multicolumn{3}{c}{\hrulefill} & \multicolumn{3}{c}{\hrulefill} &
\multicolumn{3}{c@{}}{\hrulefill}\\
$\bolds q$ & \textbf{1} & \textbf{3} & \textbf{5} & \textbf{1} & \textbf{3}
& \textbf{5} & \textbf{1} & \textbf{3} & \textbf{5} \\
\hline
$\widehat\eta_{t1}$ & 0.00 & 0.00 & 0.00 & 0.00 & 0.00 & 0.00 & 0.00 & 0.00
& 0.00 \\
$\widehat\eta_{t2}$ & 0.00 & 0.00 & 0.00 & 0.00 & 0.00 & 0.00 & 0.00 & 0.00
& 0.00 \\
$\widehat\eta_{t3}$ & 0.11 & 0.08 & 0.04 & 0.09 & 0.08 & 0.02 & 0.07 & 0.06
& 0.02 \\
$\widehat\eta_{t4}$ & 0.05 & 0.25 & 0.28 & 0.22 & 0.33 & 0.47 & 0.53 & 0.56
& 0.63 \\
$\widehat\eta_{t5}$ & 0.22 & 0.19 & 0.39 & 0.30 & 0.47 & 0.58 & 0.73 & 0.77
& 0.81 \\
\hline
\end{tabular*}
\end{table}

Figure~\ref{figure:10} displays the first $\widehat d (= 2)$ estimated
eigenfunctions $\widehat\psi_j$ in (\ref{a39}).
Although the estimated curves $Y_t(\cdot)$ in Figure~\ref{figure:8} are somehow
different for different bandwidths, the shape of the estimated
eigenfunctions is insensitive to the choice of
bandwidth.

\begin{figure}

\includegraphics{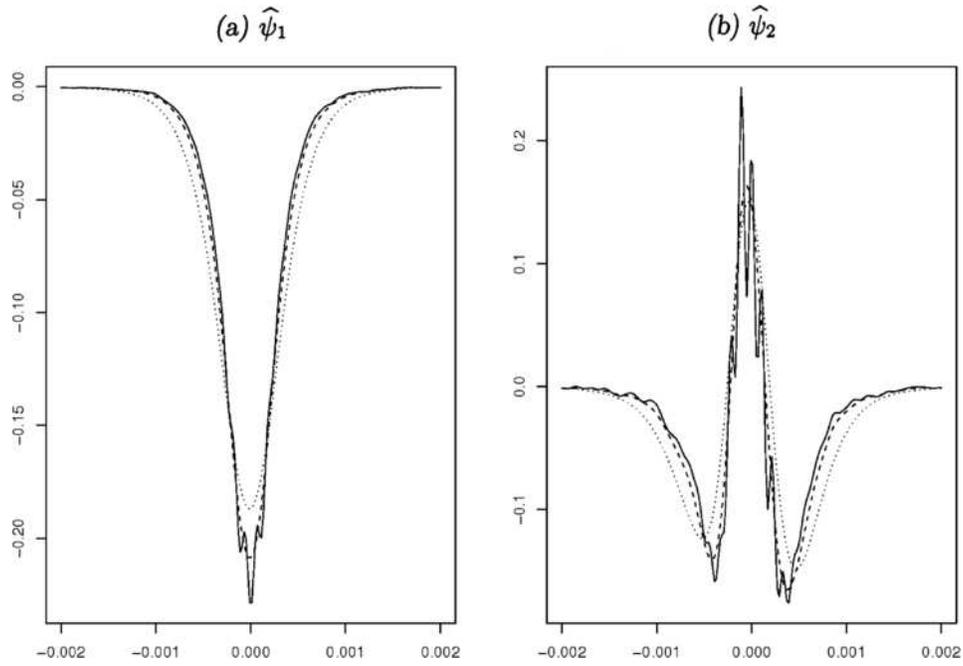}

\caption{Estimated eigenfunctions \textup{(a)}
$\widehat\psi_1$ and \textup{(b)} $\widehat\psi_2$ using bandwidths $h_t = 0.5 \widehat h_t$
(solid lines), $\widehat h_t$~(dashed lines) and $2 \widehat h_t$ (dotted
lines).}\label{figure:10}
\end{figure}

Figure~\ref{figure:11} displays time series plots of the estimated
loadings $\widehat\eta_{t1}$ and $\widehat\eta_{tj}$. Again
the estimated loadings with three levels of bandwidth
are almost indistinguishable from each other.
Furthermore, the ACF and PACF of the series $\widehat{\bfeta}_{tj}
= (\widehat\eta_{t1}, \widehat\eta_{t2})^\prime$ are also virtually identical
for all three choices of $h$. These graphics are displayed in
Figures~\ref{figure:12} and~\ref{figure:13}.

\begin{figure}

\includegraphics{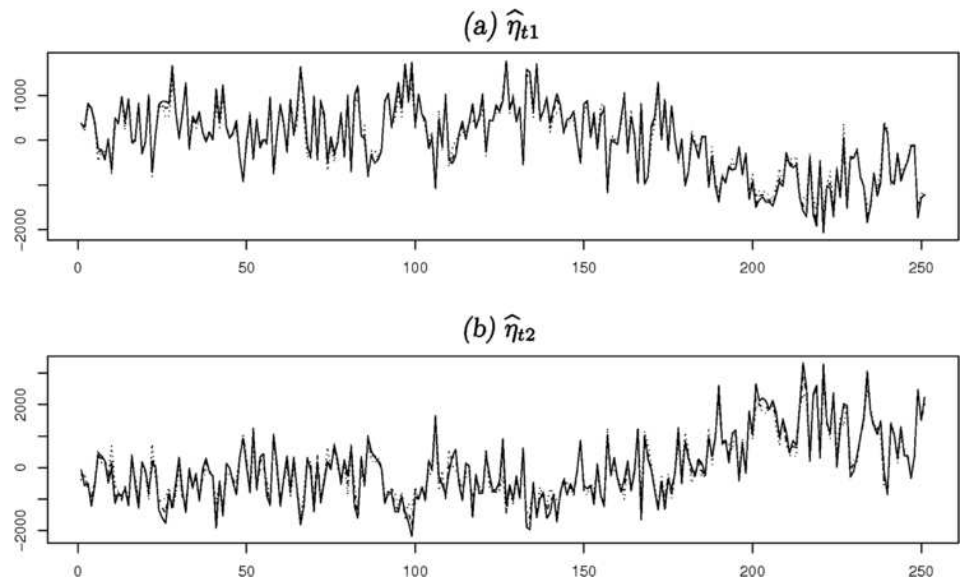}

\caption{Estimated loadings \textup{(a)} $\widehat\eta_{t1}$ and
\textup{(b)}
$\widehat\eta_{t2}$ using bandwidths $h_t = 0.5 \widehat h_t$
(solid lines), $\widehat h_t$~(dashed lines) and $2 \widehat h_t$ (dotted
lines).}\label{figure:11}
\end{figure}

\begin{figure}

\includegraphics{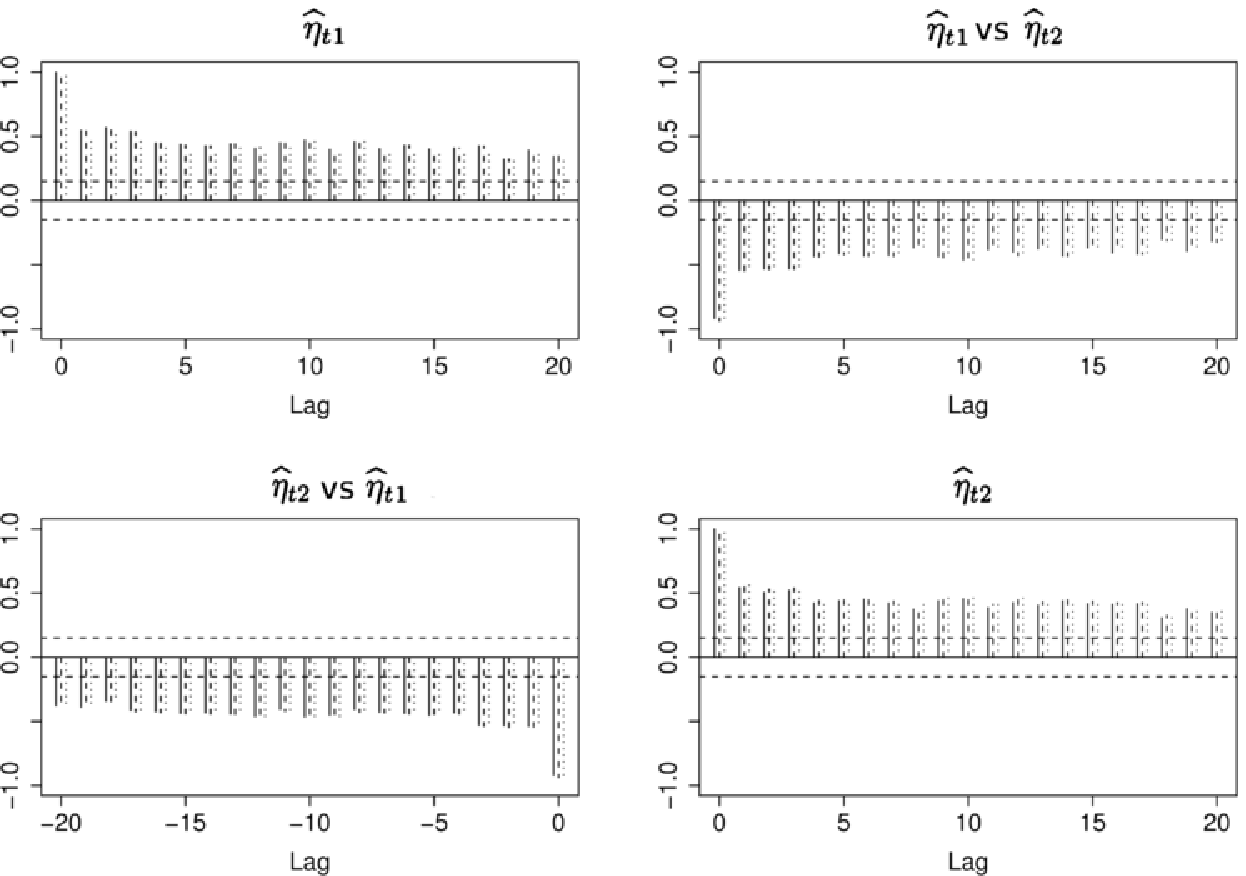}

\caption{ACF of $\widehat{\bfeta}_{tj}$ using bandwidths$h_t = 0.5 \widehat h_t$
(solid lines), $\widehat h_t$~(dashed lines) and $2 \widehat h_t$ (dotted
lines).}\label{figure:12}
\end{figure}

\begin{figure}

\includegraphics{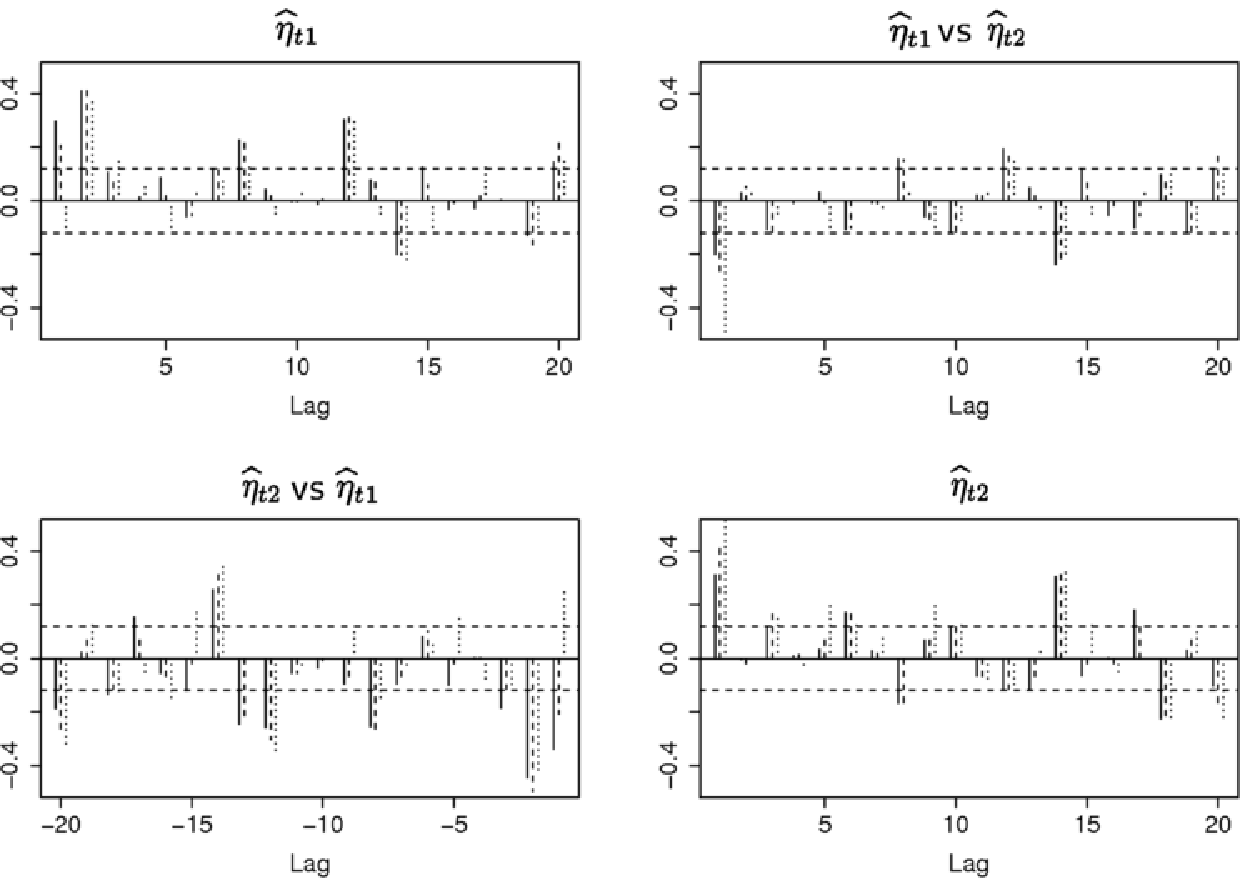}\vspace*{-3pt}

\caption{PACF of $\widehat{\bfeta}_{tj}$ using bandwidths $h_t = 0.5 \widehat h_t$
(solid lines), $\widehat h_t$ (dashed lines) and $2 \widehat h_t$ (dotted
lines).}\label{figure:13}\vspace*{-6pt}
\end{figure}

We now fit a VAR model to the estimated loadings, $\widehat{\bfeta}_{t}$:
%
\begin{equation}\label{var}
\widehat{\bfeta}_{t} = \sum_{k=1}^\tau\bA_{k} \widehat\bfeta_{t -k} + \be_t.
\end{equation}
Since the estimated loadings $\widehat\bfeta_{tj}$, as defined in (\ref{a10:2}),
have mean zero by construction, there is no intercept term in the
model. We
choose the order $\tau$ in (\ref{var}) by
minimizing the AIC. The AIC values for the order $\tau=0,1,\ldots,10$
are given in Table~\ref{table:4}. With all three bandwidths used, the AIC
chooses $\tau=3$, and
the multivariate portmanteau\vadjust{\goodbreak} test (with lag values
1, 3 and 5) of \citet{li1981} for the residual of the fitted VAR
models are insignificant at the $10\%$ level. The Yule--Walker
estimates of the parameter matrices, $\bA_k = (a_{k,ij})$ in
(\ref{var}), with the order $\tau=3$ are given in Table
\ref{table:5}.

To summarize, we found that the dynamic behavior of the IBM intraday
return densities in 2006 was driven by two factors. These
factors series are modeled well by a VAR(3) process. We note
that with all the three levels of smoothness adopted in the initial
density estimation, these conclusions were unchanged.

Finally, we make a cautionary remark on the implied true curves
$X_t(\cdot)$
in the above analysis. We take the unknown true daily densities as
$X_t(\cdot)$. We see those densities as random curves, as the
distribution of the intraday returns tomorrow depends on the distributions
of today, yesterday and so on, but is not
entirely determined by them.
Now in model (\ref{a1}), $E\{ \ve_t(u) \} = E\{ Y_t(u) \} - E\{
X_t(u) \} \not\equiv0$.
But this does not affect the analysis performed in identifying the
dimensionality of the curves; see also Pan and Yao (\citeyear{pan2008}).
Note that (\ref{a38}) provides an alternative estimator for the true density
$X_t(\cdot)$ based on the dynamic structure of the curve series. It can be
used, for example, to forecast the density for tomorrow. However,
an obvious normalization should be applied since we did not make use the
constraint $\int_\calI X_t(u) \,du =1$ in constructing
(\ref{a38}).\vadjust{\goodbreak}

\begin{appendix}
\section{}\label{appA}


In this section, we provide the relevant background on operator
theory used in this work. More detailed accounts may be found in
\citet{dunford1988}.

Let $\calH$ be a real separable Hilbert space with respect to some
inner product $\inner{\cdot}{\cdot}$. For any $\calV\subset\calH$,
the orthogonal complement of $\calV$ is given by
\[
\calV^\bot= \{x \in\calH\dvtx\inner{x}{y}=0, \forall y \in
\calV\}.
\]
Note that $\calV^{\bot\bot} = \overline{\calV}$ where
$\overline{\calV}$ denotes the closure of $\calV$. Clearly, if
$\calV$ is finite dimensional then $\calV^{\bot\bot} = \calV$.

Let $L$ be a linear operator from $\calH$ to $\calH$. For $x\in
\calH$, denote by $Lx$ the image of $x$ under $L$. The adjoint of
$L$ is denoted by $L^*$ and satisfies
%
\begin{equation}\label{A1}
\inner{Lx}{y} = \inner{x}{L^*y},\qquad x,y\in
\calH.
\end{equation}
$L$ is said to be self adjoint if $L^* = L$ and nonnegative
definite if
\[
\inner{Lx}{x} \geq0\qquad \forall x \in\calH.
\]
The image and null space of $L$ are defined as $\image{L} = \{y \in
\calH\dvtx y=Lx, x \in\calH\}$ and $\kernel{L} = \{x \in\calH\dvtx
Lx = 0\}$, respectively. Note that $\operatorname{Ker}(L^*) =
(\operatorname{Im}(L))^\bot$, $\operatorname{Ker}(L) = (\operatorname{Im}(L^*))^\bot$ and
$\operatorname{Ker}(L^*) = \operatorname{Ker}(LL^*)$. We define the rank of $L$ to
be $r(L) = \dim(\image{L})$ and we say that $L$ is finite
dimensional if $r(L) < \infty$.\vadjust{\goodbreak}

\begin{table}
\def\arraystretch{0.9}
\caption{AIC values from fitting the VAR model
in (\protect\ref{var}). The figures in this table have been centered at the
minimum AIC value}\label{table:4}
\begin{tabular*}{\tablewidth}{@{\extracolsep{\fill}}lcccccc@{}}
\hline
& $\bolds{\tau=0}$ & $\bolds{\tau=1}$ & $\bolds{\tau=2}$ & $\bolds{\tau=3}$ & $\bolds{\tau=4}$ & $\bolds{\tau=5}$\\
\hline
$h_t=0.5 \widehat h_t$ & 131.33& 40.39 & 9.98 & 0.00 & 7.86 & 10.38\\
$h_t = \widehat h_t$ & 133.04& 41.32 & 9.53 & 0.00 & 7.47 & 10.08\\
$h_t= 2 \widehat h_t$ & 135.47& 40.83 & 9.58 & 0.00 & 7.00 & \phantom{0}8.94 \\
\hline
\end{tabular*}
\end{table}

\begin{table}[b]
\def\arraystretch{0.9}
\caption{Estimated parameter matrices $\bA_k
= (a_{k,ij})$ from fitting the VAR model in (\protect\ref{var})}\label{table:5}
\begin{tabular*}{\tablewidth}{@{\extracolsep{\fill}}ld{2.2}d{2.2}d{2.2}d{2.2}d{2.2}d{2.2}@{}}
\hline
\multicolumn{1}{@{}l}{$\bolds{j}$} & \multicolumn{3}{c}{\textbf{1}} & \multicolumn{3}{c}{\textbf{2}}
\\[-4pt]
& \multicolumn{3}{c}{\hrulefill} &
\multicolumn{3}{c@{}}{\hrulefill}\\
\multicolumn{1}{@{}l}{$\bolds{h_t}$} & \multicolumn{1}{c}{$\bolds{0.5\widehat h_t}$}
& \multicolumn{1}{c}{$\bolds{\widehat h_t}$} & \multicolumn{1}{c}{$\bolds{2 \widehat
h_t}$} & \multicolumn{1}{c}{$\bolds{0.5 \widehat h_t}$}
& \multicolumn{1}{c}{$\bolds{\widehat h_t}$} & \multicolumn{1}{c@{}}{$\bolds{2 \widehat h_t}$} \\
\hline
$a_{1,1j}$ & 0.08 & 0.07 & 0.01 & -0.14 & -0.16 & -0.22 \\
$a_{1,2j}$ & -0.08 & -0.05 & 0.03 & 0.24 & 0.26 & 0.33 \\
[2pt]
$a_{2,1j}$ & 0.35 & 0.39 & 0.38 & 0.06 & 0.09 & 0.08 \\
$a_{2,2j}$ & -0.36 & -0.43 & -0.43 & -0.05 & -0.10 & -0.11 \\
[2pt]
$a_{3,1j}$ & 0.08 & 0.05 & 0.02 & -0.13 & -0.15 & -0.18 \\
$a_{3,2j}$ & -0.16 & -0.13 & -0.11 & 0.14 & 0.15 & 0.17 \\
\hline
\end{tabular*}
\end{table}

A linear operator $L$ is said to be bounded if there exists some
finite constant $\Delta>0$ such that for all $x \in\calH$
\[
\vectornorm{L x} < \Delta\vectornorm{x},
\]
where \mbox{$\vectornorm{\cdot}$} is the norm induced on $\calH$ by
$\inner{\cdot}{\cdot}$. We denote the space of bounded linear
operators from $\calH$ to $\calH$ by $\calB= \calB(\calH, \calH)$
and the uniform topology on $\calB$ is defined by
\[
\uniform{L} = \sup_{\vectornorm{x} \leq1} \vectornorm{ L
x},\qquad L \in\calB.
\]
Note that all bounded linear operators are continuous, and the
converse also holds.

An operator $L \in\calB$ is said to be compact if there exists two
orthonormal sequences $\{e_j\}$ in $\{f_j\}$ of $\calH$ and a
sequence of scalars $\{\lambda_j\}$ decreasing to zero such that
\[
Lx = \sum_{j=1}^\infty\lambda_j \inner{e_j}{x} f_j,\qquad
x\in\calH,
\]
or more compactly
%
\begin{equation}\label{A2}
L = \sum_{j=1}^\infty\lambda_j e_j \otimes f_j.\vadjust{\goodbreak}
\end{equation}
Note that if $\calH= \calL_2(\calI)$ equipped with the inner
product defined in (\ref{a36}),~then
\[
(L x) (u) = \sum_{j=1}^\infty\lambda_j \inner{e_j}{x}
f_j(u).
\]
Clearly, $\image{L} = \operatorname{sp}\{f_j \dvtx j \geq1\}$ and
$\kernel{L} = \operatorname{sp}\{e_j \dvtx j \geq1 \}^\bot$.

The Hilbert--Schmidt norm of a compact linear operator $L$ is defined
as $\hilbert{L} = (\sum_{j=1}^\infty\lambda_j^2)^{1/2}$. We
will\vspace*{2pt}
let $\calS$ denote the space consisting of all the operators with a
finite Hilbert--Schmidt or nuclear norm. Clearly, we have the
inequalities $\hilbert{\cdot} \geq\uniform{\cdot}$, and thus the
inclusions $\calS\subset\calB$. Note that $\calB$ is a Banach
space when equipped with their respective norms. Furthermore, $\calS$
is a Hilbert space with respect to the inner product
\[
\inner{L_1}{L_2}_\calS= \sum_{i,j=1}^\infty\inner{L_1g_i}{h_j}
\inner{L_2g_i}{h_j},\qquad L_1,L_2 \in\calS,
\]
where $\{g_i\}$ and $\{h_j\}$ are any orthonormal bases of $\calH$.


\section{}\label{appB}


In this section, we provide the proofs for the propositions in
Section~\ref{section:2.2} and the theorems in Section
\ref{section:2.3}. Throughout the proofs, we may use $C$ to denote
some (generic) positive and finite constant which may vary from line
to line. We introduce some technical lemmas first.
\begin{lemma}\label{lemma:1}
Let $L$ be a finite-dimensional operator such that for some
sequences of orthonormal vectors $\{e_j\}$, $\{f_j\}$, $\{g_j\}$ and
$\{h_j\}$ and some sequences of decreasing scalars $\{\theta_j\}$
and $\{\lambda_j\}$, $L$ admits the spectral decompositions $L =
\sum_{j=1}^d \theta_j e_j \otimes f_j = \sum_{j=1}^{d^\prime}
\lambda_j g_j \otimes h_j$. Then it holds that $d^\prime= d$.
\end{lemma}
\begin{pf}
Note that if $d \neq d^\prime$
then both $\image{L}$ and $\image{L_k^*}$ will be of different dimensions
under the alternative characterizations due to linear independence
of $\{e_j\}$, $\{f_j\}$, $\{g_j\}$ and $\{h_j\}$. Thus, it must hold
that $d=d^\prime$.
\end{pf}
\begin{lemma}\label{lemma:2}
Let $L$ be a linear operator from $\calH$ to $\calH$, where $\calH$
is a separable Hilbert space. Then it holds that
$\overline{\image{LL^*}} = \overline{\image{L}}$.
\end{lemma}
\begin{pf}
Using the facts about
inner product spaces and linear operators stated in Appendix~\ref{appA}, we
have
\begin{eqnarray*}
\overline{\image{LL^* }} &=& (\image{LL^*})^{\bot\bot} =
(\image{(LL^*)^*})^{\bot\bot}\\
&=&
(\kernel{LL^*})^\bot
=(\kernel{L^*})^\bot\\
&=&
(\image{L})^{\bot\bot} =
\overline{\image{L}},
\end{eqnarray*}
which concludes the proof.\vadjust{\goodbreak}
\end{pf}

For the sake of the simplicity in presentation of the proofs, we
adopt the the standard notation for Hilbert spaces. For any $f \in
\calL_2(\calI)$, we write $\vectornorm{f} = \sqrt{\inner{f}{f}}$
[see (\ref{a36})], and denote $M_k f \in\calL_2(\calI)$ the image
of $f$ under the operator $M_k$ in the sense that
\[
(M_kf)(u) = \int_\calI M_k(u,v) f(v) \,dv.
\]
The operators $N_k, K, \widehat M_k $ and $\widehat K$ may be expressed in the
same manner. Note now that the adjoint operator of $M_k$ is
\[
(M_k^*f)(u) = \int_\calI M_k(v,u) f(v) \,dv.
\]
See (\ref{A1}). Furthermore, $N_k = M_k M_k^*$ in the sense that $N_k
f= M_k M_k^*f$; see~(\ref{a34}). Similarly, $\widehat K = \sum_{k=1}^p \widehat
M_k \widehat M_k^*$; see (\ref{a35}).
\begin{pf*}{Proof of Proposition~\ref{prop:1}} (i) To save notational burden, we
set $k \equiv k_0$. We only need to show
$\image{N_k} = \calM$. Since $N_k = M_k M_k^*$, it follows from
Lemma~\ref{lemma:2} that $\image{N_k} = \image{M_k M_k^*} =
\image{M_k}$ as $N_k$ and $M_k$ are finite dimensional and thus
their images are closed.

Now, recall from Section~\ref{section:2.2.1} that $M_k$ may be
decomposed as
%
\begin{equation}\label{M*}
M_k = \sum_{i,j=1}^d \sigma_{ij}^{(k)} \varphi_i \otimes\varphi_j.
\end{equation}
See also (\ref{A2}). Thus, from (\ref{M*}), we may write
%
\begin{equation}\label{lemma1-1}
M_{k} = \sum_{i=1}^d \lambda_{i}^{(k)} \varphi_i \otimes
\rho_{i}^{(k)},
\end{equation}
where
\[
\rho_{ik}= \frac{\sum_{j=1}^d \sigma_{ij}^{(k)}
\varphi_j}{\vectornorm{\sum_{j=1}^d \sigma_{ij}^{(k)}
\varphi_j}},\qquad
\lambda_{i}^{(k)} = \Biggl\|\sum_{j=1}^d
\sigma_{ij}^{(k)} \varphi_i\Biggr\|.
\]
From (\ref{M*}), it is clear that $\image{M_k} \subseteq\calM$,
which is finite dimensional. Thus, $M_k$ is compact and therefore
admits a spectral decomposition of the form
%
\begin{equation}\label{lemma1-2}
M_k = \sum_{j=1}^{d_k} \theta_{j}^{(k)} \psi_{j}^{(k)} \otimes
\phi_j^{(k)}
\end{equation}
with $(\phi_{j}^{(k)},\psi_{j}^{(k)})$ forming the adjoint pair of
singular functions of $M_k$ corresponding to the singular value
$\theta_{j}^{(k)}$. Clearly, $d_k \leq d$. Thus, if $d_k < d$,
$\image{M_k} \subset\calM$ since from (\ref{lemma1-2}),
$\image{M_k} = \operatorname{span}\{\phi_{j}^{(k)} \dvtx j=1,\ldots,d_k\}$ and
any subset of $d_k < d$ linearly independent elements in a
$d$-dimensional space can only span a proper subset of the original
space.\vadjust{\goodbreak}

Now to complete the proof, we only need to show that the set of
$\{\rho_{j}^{(k)}\}$ in (\ref{lemma1-1}) is linearly independent for
some $k$. If this can be done, then we are in a position to apply
Lemma~\ref{lemma:1}. Let $\bbeta$ be an arbitrary vector in
$\mathbb{R}^d$ and put $\bvarphi=
(\varphi_1,\ldots,\varphi_d)^\prime$ and $\brho_k =
(\rho_{1}^{(k)},\ldots,\rho_{d}^{(k)})^\prime$, then the linear
independence of the set $\{\rho_{i}^{(k)}\}$ can easily be seen as
the equation
\[
\bbeta\brho_k = \bbeta\bSigma_k \bvarphi= 0
\]
has a nontrivial solution if and only if $\bbeta\bSigma_k =
\mathbf{0}$. However, since $\bSigma_k$ is of full rank by
assumption, it follows that it is invertible and the only solution
is the trivial one $\bbeta= \mathbf{0}$. Thus, Lemma~\ref{lemma:1}
implies $d_k = d$ and the result follows from noting that any
linearly independent set of $d$ elements in a $d$-dimensional vector
space forms a basis for that space.

(ii) Similar to the proof of part (i) above, we only need to show
$\image{K} = \calM$. Note that for any $f \in L_2(\calI)$, $\langle M_k
M_k^* f, f\rangle = \langle M_k^* f, M_k^* f\rangle = \vectornorm{M_k^* f}^2
\geq0$, thus the composition $N_k = M_k M_k^*$ is nonnegative
definite which implies that $K$ is also nonnegative definite.
Therefore, $\operatorname{Im}(K)= \bigcup_{k=1}^p \operatorname{Im}(N_k)$. From
here, the result given in part (i) of the proposition concludes the
proof.\vspace*{-3pt}
\end{pf*}
\begin{pf*}{Proof of Proposition~\ref{proposition2}}
Let $\widehat\theta_j$ be a nonzero eigenvalue of $\bK^*$, and
$\bgamma_j = (\gamma_{1j},\ldots, \gamma_{n-p, j})'$ be the
corresponding eigenvector, that is, $ \bK^* \bgamma_j =
\bgamma_j\widehat\theta_j $. Writing this equation component by
component, we obtain that
%
\begin{equation} \label{A3}\qquad
{1 \over(n-p)^2} \sum_{i,s=1}^{n-p} \sum_{k=1}^p
\inner{Y_{t+k}-\bar Y}{Y_{s+k}-\bar Y} \inner{Y_s-\bar Y}{Y_i - \bar
Y} \gamma_{ij} = \gamma_{tj} \widehat\theta_j
\end{equation}
for $t=1,\ldots,n-p$; see (\ref{K}). For $\widehat\psi_j$ defined in
(\ref{a39}),
\begin{eqnarray*}
(\widehat K \widehat\psi_j) (u) &=& \int_\calI\widehat K(u,v) \widehat\psi_j(v) \,dv\\[-3pt]
&=& {1 \over(n-p)^2} \sum_{t,s=1}^{n-p} \sum_{k=1}^p
\{ Y_t(u) - \bar Y(u) \} \inner{Y_s - \bar Y}{\widehat\psi_j}\\[-4pt]
&&\hspace*{78.1pt}{}\times \inner
{Y_{t+k}-\bar Y}{Y_{s+k}-\bar Y}\\[-3pt]
&=& {1 \over(n-p)^2} \sum_{t,s, i=1}^{n-p} \sum_{k=1}^p \{ Y_t(u) -
\bar Y(u) \} \gamma_{i j} \inner{Y_s - \bar Y}{Y_i - \bar Y}\\[-4pt]
&&\hspace*{83.3pt}{}\times
\inner{Y_{t+k}-\bar Y}{Y_{s+k}-\bar Y};
\end{eqnarray*}
see (\ref{a35}). Plugging (\ref{A3}) into the right-hand side of the
above expression, we obtain that
\[
(\widehat K \widehat\psi_j) (u) = \sum_{t=1}^{n-p} \{ Y_t(u) - \bar Y(u) \}
\gamma_{tj} \widehat\theta_j = \widehat\psi_j(u) \widehat\theta_j,
\]
that is, $\widehat\psi_j$ is an eigenfunction of $\widehat K$ corresponding to
the eigenvalue $\widehat\theta_j$.\vadjust{\goodbreak}
\end{pf*}

As we shall see, the operator $\widehat K = \sum_{k=1}^p\widehat M_k \widehat M_k^*$
may be written as a functional of empirical distributions of Hilbertian
random variables. Thus, we
require an auxiliary result to deal with this form of process. To
this end, we extend the $V$-statistic results of \citet{sen1972} to
the setting of Hilbertian valued random variables. Further details
about $V$-statistics may be found in \citet{lee1990}.

Let $\calH$ be a real separable Hilbert space with norm
$\norm{\cdot}$ generated by an inner product \mbox{$\inner{\cdot}{\cdot}$}.
Let $X_t \in\calX$ be a sequence of strictly stationary and Hilbertian
random variables whose distribution
functions will be denoted by $P(x), x \in\calH$. Note that the
spaces $\calX$ and $\calH$ may differ. Let $\phi\dvtx \calX^m
\to\calH$ be Bochner integrable and symmetric in each
of its $m$($\geq$2) arguments. Now consider the functional
\[
\theta(P) = \int_{\calX^m} \phi(x_1,\ldots,x_m) \prod_{j=1}^m
P(dx_j),
\]
defined over $\calP= \{P\dvtx \norm{\theta(P)} < \infty\}$. As an
estimator of $\theta(P)$, consider the $V$-statistic defined by
\[
V_n = n^{-m} \sum_{i_1 = 1}^n \cdots\sum_{i_m=1}^n
\phi(X_{i_1},\ldots, X_{i_m}).
\]

Now for $c=0,1,\ldots,m$, we define the functions
\[
\phi_c(x_1,\ldots,x_c) = \int_{\calX^{m-c}}
\phi(x_1,\ldots,x_c,x_{c+1},\ldots,x_m) \prod_{j=c+1}^m
P(dx_{j})
\]
and
\[
g_c(x_1,\ldots,x_c) = \sum_{d=0}^c (-1)^{c-d} \sum_{1 \leq j_1 <
\cdots< j_d \leq c} \phi_d(X_{j_1},\ldots, X_{j_d}).
\]
In order to construct the canonical decomposition
of $V_n$, we use Dirac's $\delta$-measure to define the
empirical measure $P_n$ as follows:
\[
P_n(A) = n^{-1} \bigl( \delta_{X_1}(A)+ \cdots+
\delta_{X_n}(A) \bigr),\qquad A \in\calX.
\]
Then for $c=1,\ldots,m$, we set
\begin{eqnarray*}
V_{nc} &=& \int_{X^c} \phi_c(x_1,\ldots,x_c) \prod_{j=1}^c \bigl(
P_n(dx_j) - P(dx_j) \bigr)\\
&=& n^{-c} \sum_{i_1=1}^n \cdots\sum_{i_c=1}^n
g_c(X_{i_1},\ldots,X_{i_c}),
\end{eqnarray*}
then we have
%
\begin{equation}\label{hoef-v}
V_n - \theta(P) = \sum_{c=1}^m \pmatrix{m \cr c} V_{nc}.\vadjust{\goodbreak}
\end{equation}
In particular, note that
\[
V_{n1} = \frac{1}{n} \sum_{i=1}^n g_1(X_i).
\]
Decomposition \eqref{hoef-v} is the Hoeffding representation of the
statistic $V_n$. It plays a central
role in the proof of Lemma~\ref{lemma:3} below. We are now in a
position to state some regularity conditions which form the basis of
the result.

\begin{enumerate}[$\bullet$ A3.]
\item[$\bullet$ A1.] $\{X_t\}$ is strictly stationary and $\psi$-mixing with
$\psi$-mixing coefficients satisfying the condition $\sum_{l=1}^\infty
l^{m-1} \psi^{1/2}(l) < \infty$.\vspace*{2pt}
\item[$\bullet$ A2.] $\int_{\calX^m} \norm{\phi(x_1,\ldots,x_m)}^{2}
\prod_{j=1}^m P(dx_j) < \infty.\nonumber$
\item[$\bullet$ A3.] $E\norm{g_1 (X_1)}^2 + 2 \sum_{k=2}^\infty E\inner
{g_1(X_1)}{g_1{X_k}} \neq
0$.
\end{enumerate}
\begin{lemma}\label{lemma:3}
Let conditions \textup{A1--A3} hold. Then for $c=1,\ldots,m$ it holds that $E\norm
{V_{nc}}^2 = O(n^{-c})$.
\end{lemma}
\begin{pf}
We make use of
(\ref{hoef-v}). Let $\{e_j \dvtx j \geq
1\}$ be an orthonormal basis of $\calH$. Then
%
\begin{equation}\label{prop3:1}
E \norm{ V_{nc}}^2 = \sum_{j=1}^\infty E \inner{e_j}{V_{nc}}^2,
\end{equation}
where $\inner{e_j}{V_{nc}}$ is the $\mathbb{R}$ valued $V$-statistic
\[
\inner{e_j} {V_{nc}} = n^{-c} \sum_{i_1=1}^n \cdots\sum_{i_c=1}^n
\inner{e_j}{g_c(X_{i_1},\ldots,X_{i_c})}.
\]

Now under conditions A1--A3, Lemma 3.3 in \citet{sen1972} yields
%
\begin{equation}\label{prop3:2}
E \inner{e_j} {V_{nc}}^2 \leq C n^{-c} \int_{\calX^c} \inner{e_j}
{\phi_c
(x_1,\ldots,x_c)}^{2} \prod_{j=1}^c P(dx_j)
\end{equation}
for all $j \geq1$. Now inserting the estimate in (\ref{prop3:2})
into (\ref{prop3:1}) yields
\begin{eqnarray*}
E \norm{ V_{nc}}^2 &\leq& C n^{-c} \sum_{j=1}^\infty
\int_{\calX^c} \inner{e_j} {\phi_c
(x_1,\ldots,x_c)}^{2} \prod_{j=1}^c
P(dx_j)\\
&\leq& C n^{-c} \int_{\calX^c} \norm{\phi_c
(x_1,\ldots,x_c)}^{2} \prod_{j=1}^c
P(dx_j)\\
&\leq& C n^{-c} \sum_{j=1}^\infty\int_{\calX^c} \norm{\phi
(x_1,\ldots,x_m)}^{2} \prod_{j=1}^m
P(dx_j)\\
&=& O(n^{-c})
\end{eqnarray*}
as required.\vadjust{\goodbreak}
\end{pf}
\begin{pf*}{Proof of Theorem~\ref{theorem1}}
(i) Since $p$ is fixed and
finite, we may set $n \equiv n-p$. Let $Z_{tk} = (Y_t - \mu) \otimes
(Y_{t+k} - \mu) \in
\calS$. Now consider the kernel $\rho\dvtx \calS\times\calS
\to\calS$ given by
%
\begin{equation}\label{u:kernel}
\rho(A, B) = A B^*,\qquad A, B \in\calS.
\end{equation}
Now note that from (\ref{u:kernel}),
\[
\widehat M_k \widehat M_k = n^{-2} \sum_{i=1}^n \sum_{j=1}^n \rho( Z_{ik},
Z_{jk}),
\]
which in light of the preceding discussion is simply a $\calS$
valued von Mises functional. Then $d \geq1$ it holds that $M_k \neq
0$, an application of Lemma~\ref{lemma:3} yields
%
\begin{equation}\label{th1:1}
E\norm{\widehat M_k \widehat M_k^* - M_k M_k^*}_\calS^2 = O (n^{-1}).
\end{equation}
Note that if $d = 0$, the rate in (\ref{th1:1}) would be $n^{-2}$, that
is, the kernel $\rho$ would possess the property of first order
degeneracy. Now by (\ref{th1:1}) and the Chebyshev inequality, we have
\[
\norm{\widehat K - K}_\calS\leq\sum_{k=1}^p \norm{\widehat M_k \widehat M_k - M_k
M_k^*}_\calS= O_p(n^{-1/2}).
\]

\mbox{\phantom{i}}(ii) Given $\hilbert{\widehat K - K} = O_p(n^{-1/2})$, Lemma 4.2 in
\citet{Bosq2000} implies the $\sup_{j \geq1} | \widehat\theta_j -
\theta_j| \leq\norm{\widehat K - K}_\calS= O_p(n^{-1/2})$. Condition
C3
ensures that $\psi_j$ is an
identifiable statistical parameter for $j=1,\ldots,d$. From here,
Lemma 4.3 in \citet{Bosq2000} implies $\norm{\widehat\psi_j - \psi_j}
\leq C \norm{\widehat K - K}_\calS= O_p(n^{-1/2})$.

(iii) First, note that by Lemma~\ref{lemma:3} we have
%
\begin{equation}\label{th1:2.1}
E \hilbert{ \widehat M_k \widehat M_k^* - \widehat M_k M_k^*}^2 = O(n^{-2}).
\end{equation}
Put $\widetilde K = \sum_{k=1}^p \widehat M_k M_k$. Then by (\ref{th1:2.1}) and the
Chebyshev inequality, we have
%
\begin{equation}\label{th1:2}
\norm{ \widehat K - \widetilde K}_\calS\leq\sum_{k=1}^p \norm{\widehat M_k \widehat
M_k^* - \widehat M_k M_k^*}_\calS= O_p(n^{-1}).
\end{equation}
The estimate in (\ref{th1:2}) will prove to be crucial in deriving the
results for $\widehat\theta_j$ when $j \geq d+1$.

Now, extend $\psi_1, \ldots, \psi_d$ to a complete orthonormal basis
of $\calH$. Then it holds that
%
\begin{equation}\label{th1:4}
\sum_{j=1}^n \widehat\theta_j = \sum_{j=1}^\infty\inner{\psi_j}{\widehat K
\psi_j},
\end{equation}
and by recalling that $\theta_j = 0$ for $j >d$
%
\begin{equation}\label{th1:5}
\sum_{j=1}^d \theta_j = \sum_{j=1}^d \inner{\psi_j} {K \psi_j}.\vadjust{\goodbreak}
\end{equation}
Note that $\operatorname{span} \{\psi_j \dvtx j > d\} = \calM^\bot$ and $K
\psi_j = 0$ for all $j >d$ since $\kernel{K} = \calM^\bot$. Thus,
from (\ref{th1:4}) and (\ref{th1:5}), we have
%
\begin{equation}\label{th1:6}
\sum_{j=1}^n \widehat\theta_j - \theta_j = \sum_{j=1}^\infty
\inner{\psi_j}{(\widehat K - K)\psi_j}.
\end{equation}

Now we will show that
%
\begin{equation}\label{th1:7}
\widehat\theta_j - \theta_j = \inner{\psi_j}{(\widehat K - K) \psi_j} +
O_p(n^{-1}),\qquad j=1,\ldots,d.
\end{equation}
Let $K_j = \inner{\psi_j}{(\widehat K - K) \widehat\psi_j}$. Then using the
relations $K \psi_j = \theta_j\psi_j$ and $\widehat K \widehat\psi_j =
\widehat\theta_j\widehat\psi_j$ along with the fact that $K$ is self adjoint,
we have
%
\begin{eqnarray}\label{th1:8}
| K_j - (\widehat\theta_j - \theta_j) | &=& | \inner{\psi_j}{\widehat K
\widehat\psi_j} - \inner{K \psi_j}{\widehat\psi_j} - (\widehat\theta_j
-\theta_j)|\nonumber\\
&=& | (\widehat\theta_j - \theta_j) (\inner{\psi_j}{\widehat\psi_j} - 1)|
\\
&=& |\widehat\theta_j - \theta_j| | \inner{\psi_j}{\widehat\psi_j} - 1|.\nonumber
\end{eqnarray}
Note that
%
\begin{equation}\label{th1:9}\qquad
| \inner{\psi_j}{\widehat\psi_j} - 1| = | \inner{\psi_j}{\widehat\psi_j -
\psi_j}| \leq\norm{\psi_j}\norm{\widehat\psi_j - \psi_j} =
\norm{\widehat\psi_j - \psi_j}.
\end{equation}
Thus, from the results in (b) above (\ref{th1:8}) and (\ref{th1:9}),
we have $| K_j - (\widehat\theta_j - \theta_j) | \leq|\widehat\theta_j -
\theta_j| \norm{\widehat\psi_j - \psi_j} = O_p(n^{-1})$ for
$j=1,\ldots,d$.

Next, we have
\begin{eqnarray*}
| \inner{\psi_j}{(\widehat K - K)\psi_j} - K_j| &=& | \inner{ \psi_j -
\widehat\psi_j}{(\widehat K - K)\psi_j}|\\
&\leq& \norm{\psi_j - \widehat\psi_j}\norm{(\widehat K - K) \psi_j}\\
&\leq& \|\psi_j - \widehat\psi_j\|\|\widehat K - K\|_\calS,
\end{eqnarray*}
from which the results in (i) and (ii) $| \inner{\psi_j}{(\widehat K -
K)\psi_j} - K_j| = O_p(n^{-1})$, thus proving (\ref{th1:7}).

Now from (\ref{th1:7}) we have
\[
\sum_{j=1}^d \widehat\theta_j - \theta_j = \sum_{j=1}^d
\inner{\psi_j}{(\widehat K - K)\psi_j} + O_p(n^{-1}),
\]
and thus from (\ref{th1:2}) and (\ref{th1:6})
\begin{eqnarray*}
\sum_{j=d+1}^n \widehat\theta_j &=& \sum_{j=d+1}^\infty
\inner{\psi_j}{(\widehat K - K)\psi_j} + O_p(n^{-1}) \\
&=&
\sum_{j=d+1}^\infty\inner{\psi_j}{(\widetilde K - K)\psi_j} +
O_p(n^{-1}).
\end{eqnarray*}
By noting that $\psi_j \in\calM^\bot$ for $j \geq d+1$ and
$\kernel{M_k} = \kernel{\widetilde K} = \kernel{K} = \calM^\bot$, it holds
that $\sum_{j=d+1}^\infty\inner{\psi_j}{(\widetilde K - K) \psi_j} = 0$.
Thus, $\sum_{j=d+1}^n \widehat\theta_j = O_p(n^{-1})$ and the result
follows from noting that $\widehat\theta_i \leq\sum_{j=d+1}^n \widehat
\theta_j$ for $i=1,\ldots,d$.

(iv) Let $\Pi_\calM$ and $\Pi_{\calM^\bot}$ denote the projection
operators onto $\calM$ and $\calM^\bot$, respectively. Since $x =
\Pi_\calM(x) + \Pi_{\calM^\bot} (x) $ for any $x \in
\calL_2(\calI)$, we have
%
\begin{equation}\label{th1:10}
\norm{ \Pi_\calM(\widehat\psi_{i})}^2 = \norm{ \widehat\psi_{i} -
\Pi_{\calM^\bot} (\widehat\psi_{i})}^2 = \sum_{j=1}^d \inner{\widehat
\psi_{i}} {\psi_j}^2
\end{equation}
for all $i \geq1$. Now note that for $i \geq d+1$
%
\begin{eqnarray}\label{th1:11}
\norm{K(\widehat\psi_{i}) } &=& \norm{ (K - \widehat K)(\widehat\psi_{i})
+ \widehat\psi_{i} \widehat\theta_{i}}\nonumber\\
&\leq& \norm{ (K - \widehat K) (\widehat\psi_{i}) } + | \widehat
\theta_{i}|\norm{\widehat\psi_{i}}\\
&\leq& 2 \norm{ K - \widehat K}_\calB,\nonumber
\end{eqnarray}
where the final inequality follows from the definition of $\norm{
\cdot}_\calB$ and Lemma~4.2 in \citet{Bosq2000} by noting that
$\theta_{i} = 0$ for all $i \geq d+1$.

Next, we have for $i \geq d+1$
%
\begin{eqnarray}\label{th1:12}
\norm{K(\widehat\psi_{i})}^2 &=& \sum_{j=1}^\infty\inner{ K (\widehat
\psi_{i}) }{\psi_j}^2
= \sum_{j=1}^\infty\theta_j ^2 \inner{ \widehat\psi_{i}}
{\psi_j}^2\nonumber\\[-8pt]\\[-8pt]
&=& \sum_{j=1}^d \theta_j^2 \inner{ \widehat\psi_{i}}
{\psi_j}^2
\geq\theta_d^2 \sum_{j=1}^d \inner{ \widehat\psi_{i}}
{\psi_j}^2,\nonumber
\end{eqnarray}
since $\theta_1 > \cdots> \theta_d$. Combining \eqref{th1:10}, \eqref
{th1:11} and \eqref{th1:12} yields
\[
\norm{ \Pi_\calM(\widehat\psi_{d_0 +1})}^2 = \norm{ \widehat\psi_{d_0 +1}
- \Pi_{\calM^\bot} (\widehat\psi_{d_0 +1})}^2 \leq C \norm{ K - \widehat
K}_\calB,
\]
from which (i) yields the result.
\end{pf*}
\begin{lemma}\label{lemma4}
The function $D$ defined in (\ref{metric}) is a well-defined
distance measure on $\calZ_D$.
\end{lemma}
\begin{pf}
Nonnegativity, symmetry and the
identity of indiscernibles are obvious. It only remains to prove the
subadditivity property. For any $L \in\calS$, note that
$\hilbert{L} = \sqrt{\operatorname{tr}(L^*L)}$, where $\operatorname{tr}$ denotes
the trace operator. Now, for any $\calX_i \in\calZ$, $i=1,2,3$, let
$\Pi_{\calX_i}$ denote its corresponding $d$-dimensional projection
operators defined as follows:
\[
\Pi_{\calX_i} = \sum_{j=1}^d \zeta_{ij} \otimes\zeta_{ij},
\]
where $\{\zeta_{ij} \dvtx j=1,\ldots,d\}$ is some orthonormal basis of
$\calX_i$. Now the triangle inequality for the Hilbert--Schmidt norm yields
\[
\hilbert{\Pi_{\calX_1} - \Pi_{\calX_3}} \leq\hilbert{\Pi_{\calX_1}
- \Pi_{\calX_2}} + \hilbert{\Pi_{\calX_2} - \Pi_{\calX_3}}.
\]
Since the projection operators are self adjoint, we have
\begin{eqnarray*}
&&\sqrt{\operatorname{tr}(\Pi_{\calX_1}^2) + \operatorname{tr}(\Pi_{\calX_3}^2) -
2\operatorname{tr}(\Pi_{\calX_1}\Pi_{\calX_3})}\\
&&\qquad \leq \sqrt{\operatorname{tr}(\Pi_{\calX_1}^2) +
\operatorname{tr}(\Pi_{\calX_2}^2) -
2\operatorname{tr}(\Pi_{\calX_1}\Pi_{\calX_2})}\\
&&\qquad\quad{}+ \sqrt{\operatorname{tr}(\Pi_{\calX_2}^2) + \operatorname{tr}(\Pi_{\calX_3}^2) -
2\operatorname{tr}(\Pi_{\calX_2}\Pi_{\calX_3})}.
\end{eqnarray*}
Now $\operatorname{tr}(\Pi_{\calX_i}^2) = \operatorname{tr}(\Pi_{\calX_i}) = d$ and
$\operatorname{tr}(\Pi_{\calX_i}\Pi_{\calX_j}) = \sum_{k,l=1}^d \inner{\zeta
_{ik}}{\zeta_{jl}}^2$ for $i,j =1,2,3$.
These last facts along with the definition of $D$ in (\ref{metric})
give
\[
D(\calX_1, \calX_3) \leq D(\calX_1,\calX_2) +
D(\calX_2,\calX_3),
\]
which concludes the proof.
\end{pf}
\begin{pf*}{Proof of Theorem~\ref{theorem2}}
From the definition of $D$ in (\ref{metric}), note that
%
\begin{equation}\label{th2-1}
\sqrt{2d} D(\widehat\calM, \calM) = \Vert\Pi_{\widehat\calM} -
\Pi_{\calM}\Vert_\calS,
\end{equation}
where $\Pi_{\widehat\calM} = \sum_{j=1}^d \widehat\psi_j \otimes\widehat\psi_j$
and $\Pi_{\calM} = \sum_{j=1}^d\phi_j \otimes\phi_j$ with
$\phi_1,\ldots,\phi_d$ forming any orthonormal basis of $\calM$. Now
if $\Pi_\calM^1$ and $\Pi_\calM^2$ are any projection operators onto
$\calM$, then by virtue of Lemma~\ref{lemma4} it holds that
$\Vert\Pi_\calM^1 - \Pi_\calM^2\Vert_\calS= \sqrt{2d} D(\calM, \calM) =
0$. Thus, we may proceed as if $\Pi_\calM$ in (\ref{th2-1}) was
formed with eigenfunctions of $K$, that is, $\phi_j = \psi_j$ for
$j=1,\ldots,d$.

Now, we have
%
\begin{equation}\label{th2-2}\quad
\Biggl\Vert\sum_{j=1}^d \widehat\psi_j \otimes\widehat\psi_j - \sum_{j=1}^d \psi_j
\otimes\psi_j\Biggr\Vert_\calS\leq\sum_{j=1}^d \Vert\widehat\psi_j \otimes\widehat
\psi_j - \psi_j \otimes\psi_j\Vert_\calS,
\end{equation}
that is, $\widehat\psi_j \otimes\widehat\psi_j$ (resp., $\psi_j \otimes\psi_j$)
is the projection operator onto the eigensubspace generated by $\widehat
\theta_j$ (resp., $\theta_j$). Now by part (i) of Theorem
\ref{theorem1}, $\hilbert{\widehat K - K} =
O_p(n^{-1/2})$. Thus, Theorem 2.2 in \citet{Mas2003} implies that
$\Vert\widehat\psi_j \otimes\widehat\psi_j - \psi_j \otimes\psi_j\Vert_\calS=
O_p(n^{-1/2})$ for $j=1,\ldots,d$. This last fact along with~(\ref{th2-1}) and (\ref{th2-2}) yield $D(\widehat\calM, \calM) =
O_p(n^{-1/2})$.
\end{pf*}
\begin{pf*}{Proof of Theorem~\ref{theorem3}} We first note that from
(\ref{th1:1}), the triangle inequality and the $c_r$ inequality, we have
%
\begin{equation}\label{th3:1}
E\hilbert{\widehat K - K}^2 = O(n^{-1}).
\end{equation}
As $\widehat\theta_1 \geq\widehat\theta_2 \geq\cdots\geq0$ (with strict
inequality holding with probability one), it holds that $\{\widehat d > d\} =
\{\widehat\theta_{d+1} > \epsilon\}$. Now since $\theta_{d+1} = 0$, it holds that
$\widehat\theta_{d+1} = |\widehat\theta_{d+1} - \theta_{d+1}| \leq\hilbert{\widehat K -
K}$ by Lemma 4.2 in \citet{Bosq2000}. Collecting these last few facts and
applying the Chebyshev inequality yields
%
\begin{equation}\label{th3:2}
P(\widehat d > d) \leq\epsilon^{-2} E\hilbert{\widehat K - K}^2 = O((\epsilon^2 n)^{-1})
\end{equation}
by \eqref{th3:1}.
Next, we turn to $P(\widehat d < d)$.
Due to the ordering of the eigenvalues, it holds that $\{\widehat d < d\} =
\{\widehat\theta_{d-1} < \epsilon\}$. Therefore,
%
\begin{eqnarray}\label{th3:3}
P(\widehat d < d) &=& P(\widehat\theta_{d-1} < \epsilon)\nonumber\\
&=& P(\theta_{d-1} - \widehat\theta_{d-1} > \theta_{d-1} -
\epsilon)\nonumber\\[-8pt]\\[-8pt]
&\leq& P(|\theta_{d-1} - \widehat\theta_{d-1}| > \theta_{d-1} - \epsilon
)\nonumber\\
&\leq& P(\hilbert{\widehat K - K} > \theta_{d-1} - \epsilon),\nonumber
\end{eqnarray}
where the final inequality follows from Lemma 4.2 in \citet{Bosq2000}.
Now since $\theta_{d-1} > 0$ and $\epsilon\to0$ as $n
\to\infty$, it holds that $\theta_{d-1} - \epsilon> 0$
for large enough $n$. Thus, by \eqref{th3:2} and an application of
the Chebyshev inequality to \eqref{th3:1}, we have
%
\begin{equation}\label{th3:4}
P(\widehat d < d) \leq(\theta_{d-1} - \epsilon)^{-2} E\hilbert{\widehat K - K}^2
= O(n^{-1}).
\end{equation}
From \eqref{th3:2} and \eqref{th3:3}, it follows that
\[
P(\widehat d \neq d) = P(\widehat d < d) + P(\widehat d > d) = O((\epsilon^2 n)^{-1})
\to0.
\]
This completes the proof.
\end{pf*}
\end{appendix}

\section*{Acknowledgments}
We are grateful to the Associate Editor and the two referees for
their helpful comments and suggestions.

\printaddresses

\end{document}